\documentclass{ametsocV6.1}

\usepackage[inline]{trackchanges} 
\usepackage{soul}
\usepackage{subfiles}
\usepackage{subcaption}
\usepackage{commath}
\usepackage{comment}
\usepackage{multirow}
\usepackage{enumitem}
\usepackage{nicefrac}
\usepackage{tikz}
\usetikzlibrary{shapes, 
                decorations.pathreplacing, 
                calligraphy, 
                positioning}

\addeditor{Jeremy}
\addeditor{Giacomo}
\addeditor{Darren}
\addeditor{Mark}
\addeditor{Bob}

\soulregister\citep7
\soulregister\ref7
\soulregister\citet7

\renewcommand{\phi}{\varphi}
\renewcommand{\O}{\mathcal{O}}
\renewcommand{\u}{\mathbf{u}}
\newcommand{\U}{\mathbf{U}}
\newcommand{\w}{\mathbf{w}}
\newcommand{\numax}{\nu^\text{max}}
\newcommand{\numaxbetas}{\numax_{(\beta_1, \beta_2, \beta_3)}}

\makeatletter
\renewcommand{\p@subsection}{\thesection.}
\renewcommand{\p@subsubsection}{\thesection.\thesubsection.}
\makeatother

\definecolor{myred}{HTML}{ee4035}
\definecolor{myorange}{HTML}{f37736}
\definecolor{myyellow}{HTML}{fdf498}
\definecolor{mygreen}{HTML}{7bc043}
\definecolor{myblue}{HTML}{0392cf}

\title{CFL Optimized Forward-Backward Runge-Kutta Schemes for the Shallow Water Equations}

\authors{
    Jeremy R. Lilly,\aff{a}\correspondingauthor{Jeremy R. Lilly, lillyj@oregonstate.edu}
    and Darren Engwirda,\aff{b}
    and Giacomo Capodaglio,\aff{c}
    and Robert L. Higdon,\aff{a}
    and Mark R. Petersen\aff{c}
}

\affiliation{
    \aff{a}{Oregon State University, Department of Mathematics, Corvallis, OR, 97331, USA} \\
    \aff{b}{Fluid Dynamics and Solid Mechanics Group, Los Alamos National Laboratory, NM, 87545, USA} \\
    \aff{c}{Computational Physics and Methods Group, Los Alamos National Laboratory, NM, 87545, USA}
}

\abstract{
    We present the formulation and optimization of a Runge-Kutta-type time-stepping scheme for solving the shallow water equations, aimed at substantially increasing the effective allowable time-step over that of comparable methods.
    This scheme, called FB-RK(3,2), uses weighted forward-backward averaging of thickness data to advance the momentum equation.
    The weights for this averaging are chosen with an optimization process that employs a von Neumann-type analysis, ensuring that the weights maximize the admittable Courant number.
    Through a simplified local truncation error analysis and numerical experiments, we show that the method is at least second order in time for any choice of weights and exhibits low dispersion and dissipation errors for well-resolved waves.
    Further, we show that an optimized FB-RK(3,2) can take time-steps up to 2.8 times as large as a popular three-stage, third-order strong stability preserving Runge-Kutta method in a quasi-linear test case.
    In fully nonlinear shallow water test cases relevant to oceanic and atmospheric flows, FB-RK(3,2) outperforms SSPRK3 in admittable time-step by factors between roughly between 1.6 and 2.2, making the scheme approximately twice as computationally efficient with little to no effect on solution quality.
}

\begin{document}


\maketitle






\statement
The purpose of this work is to develop and optimize time-stepping schemes for models relevant to oceanic and atmospheric flows.
Specifically, for the shallow water equations we optimize for schemes that can take time-steps as large as possible while retaining solution quality.
We find that our optimized schemes can take time-steps between 1.6 and 2.2 times larger than schemes that cost the same number of floating point operations, translating directly to a corresponding speedup.
Our ultimate goal is to use these schemes in climate-scale simulations.



\section{Introduction}
\label{sec:introduction}

Large scale climate models, such as the U.S. Department of Energy's Energy Exascale Earth System Model (E3SM) \citep{golaz2022}, grow in complexity every year as new numerical techniques are developed and as high-performance computing (HPC) systems become more powerful.
Complex, climate-scale simulations can take days or weeks or more of real time to complete; it is vital that these run times are controlled without making undue sacrifices to overall model quality.
In the present work, we develop a way to formulate and optimize a class of time-stepping schemes for the shallow water equations (SWEs).
Our ultimate goal is to improve the efficiency of simulations at the climate scale, but here we focus on the shallow water equations as a starting point.  
Although this work is directed to models of ocean and atmosphere circulation, the ideas and methods developed here may generalize to other coupled systems. 

The central issue here is the following.
The Courant–Friedrichs–Lewy (CFL) condition is a well known necessary condition for the stability of a given model.
In essence, the CFL condition places a restriction on the size of the model's time-step, where this restriction depends on the model problem itself and the chosen time and space discretizations.
Formally, it is necessary for stability that
\begin{linenomath}
\begin{equation}
    \nu = c\frac{\Delta t}{\Delta x} \leq \numax, \label{eqn:cfl_condition}
\end{equation}
\end{linenomath}
where \( c \) is a speed such as a gravity wave speed. 
We refer to \( \nu \) as the Courant number and \( \numax \) as the maximal Courant number for the given situation.
In this work, we develop time-stepping schemes optimized to have maximal \( \numax \).
We do this by starting with a three-stage third-order Runge-Kutta method from \citet{wicker2002}, presently used in the Model for Prediction Across Scales-Atmosphere (MPAS-A) \citep{skamarock2012}.
This scheme is augmented by adding forward-backward averaging of the thickness tendencies used to update the fluid velocity.
Resulting schemes have large \( \numax \) and good dispersion and dissipation relations, and they compare well against other schemes in both solution accuracy and effective maximal Courant number.

The SWEs serve as a standard test problem for algorithms and methods designed for applications to ocean and atmosphere modeling; in both cases, the model equations used for climate-scale simulations are generalizations of the SWEs, in which additional vertical tendencies and forcing terms are included.
In recent years, strong stability preserving Runge-Kutta (SSPRK) methods \citep{gottlieb1998} for solving the SWEs have become quite popular; they have the property that at each time level, a spatial norm of the solution is less than or equal to the norm of the solution at the preceding level.
These methods are used to solve the SWEs in a variety of numerical paradigms and applications, such as with finite element spatial discretizations \citep{azerad2017, dawson2013, fu2022}, finite volume discretizations \citep{bao2014, capodaglio2022, hoang2019, katta2015, lilly2023, roullet2022, ullrich2010}, and layered models that use a time-splitting of barotropic and baroclinic motions \citep{karna2018, lan2022}.
Despite this widespread use, strong stability preserving (SSP) time-integration methods may not be the best choice for most large-scale geophysical SWE applications.
The SSP property is useful in systems where strong shocks and discontinuities can develop and propagate, but such dynamics are rare in climate-scale oceanic and atmospheric flows, in which barotropic evolution is typically characterised by Rossby and Kelvin-type wave modes. Numerous contemporary large-scale ocean and atmosphere models employ time-stepping methods that do not have the SSP property, motivating an investigation of its optimality for this class of flow.
For example, to solve the barotropic equations in split barotropic/baroclinic systems, the Model for Prediction Across Scales-Ocean (MPAS-O) \citep{ringler2013, petersen2019} and Model for Prediction Across Scales-Atmosphere (MPAS-A) \citep{skamarock2012} use a two-stage predictor-corrector method and a three-stage third-order Runge-Kutta method respectively.
The Community Atmosphere Model (CAM) \citep{collins2004}, the Parallel Ocean Program (POP) \citep{kerbyson2005}, the Nucleus for European Modelling of the Ocean (NEMO) \citep{madec2017}, and the HYbrid Coordinate Ocean Model (HYCOM) \citep{chassignet2007} all use variations of semi-implicit leapfrog time integration schemes for their barotropic systems.
The Regional Ocean Modeling System (ROMS) \citep{shchepetkin2005} uses a multi-step forward-backward Adams-Bashforth/Adams-Moulton method.

With the potential performance benefits of the ability to take larger time-steps in mind, we turn our attention towards methods that are optimal in this sense.
Taking advantage of the structure of the SWEs, we employ a strategy similar to that used by \citet{mesinger1977}, wherein the thickness equation is advanced in time with a simple forward Euler step, then the momentum equation is advanced with an explicit backward Euler step using the recently obtained data for the thickness.
The method, referred to by Mesinger as the forward-backward method, was more computationally efficient than popular leapfrog methods of the time.
The idea of advancing the variables of a system of ODE separately, using the most recently computed data from previously advanced variables, has been employed by the developers of ROMS.
Specifically, \citet{shchepetkin2005, shchepetkin2009} present an extension of a two-stage second-order Runge-Kutta scheme that, within each RK stage, first advances the thickness variable, then uses a weighted average of the most recently computed data with previously computed data to advance the momentum equation in a shallow water system.
The coefficients of this scheme are selected to optimize the allowable time-step length for a one-dimensional gravity wave system based on a von Neumann type stability analysis.
In this work, we develop a new, forward-backward weighted three-stage Runge-Kutta method with favorable CFL limits using a similar approach.
The optimal weights for our forward-backward scheme are obtained by von Neumann-type analysis on the two-dimensional system, including a linear Coriolis term and linearization of the advection operators about a non-zero mean state.
Using a numerical optimization algorithm, we obtain a scheme that is very efficient in terms of its maximal Courant number, as well as its dissipation and dispersion characteristics.

The paper is structured as follows.
We begin by introducing our new time-stepping scheme, referred to as FB-RK(3,2), for a general system of ODEs.
Then, we describe the process by which we obtain optimal weights for the FB averaging within each RK stage.
This process includes the derivation of a von Neumann formulation of the linearized SWEs with a nonzero mean flow and the formulation of an optimization problem.
Next, we present the results of this optimization and give a brief discussion of the dispersion and dissipation errors of the optimal FB-RK(3,2) schemes.
Finally, we perform a number of numerical experiments that demonstrate the computational efficiency and efficacy of the optimal FB-RK(3,2) schemes as compared to the three-stage third-order SSPRK scheme (SSPRK3) from \citet{gottlieb1998}.
These experiments include an investigation of CFL performance across five shallow water test cases, convergence tests, and comparisons of solution quality.
We show that an optimized FB-RK(3,2) scheme can take time-steps up to 2.2 times as large as SSPRK3 in nonlinear test cases while retaining solution quality, and because both schemes consist of three RK stages, this translates directly to a corresponding increase in efficiency.


\section{CFL Optimization of FB-RK(3,2)}
\label{sec:optimization}

Here, we introduce a time-stepping scheme for solving the shallow water equations (SWEs) that uses weighted forward-backward averaging in the layer thickness data to advance the momentum equation.
We then use a von Neumann analysis and a numerical optimization scheme to tune the forward-backward weights so that the scheme has a large Courant number when applied to the linearized SWEs.


\subsection{FB-RK(3,2)}
\label{subsec:fb-rk32}

The time-stepping scheme presented here is an extension of the three-stage, third-order Runge-Kutta time-stepping scheme RK3 from \citet{wicker2002}.
When applied to the SWEs, this extension allows us to use the most recently obtained data for the layer thickness to update the momentum data within each Runge-Kutta stage of the method.
This is done by taking a weighted average of layer thickness data at the old time level \( t^n \) and the most recent RK stage, then applying this to the momentum equation.
For the remainder of the manuscript, we will call our extension of RK3 that uses forward-backward averaging FB-RK(3,2), as it is a three-stage, second-order scheme.

Consider a general system of ODEs in independent variables \( u = u(t) \) and \( h = h(t) \) of the form
\begin{linenomath}
\begin{equation}
\begin{aligned}
    \od{u}{t} &= \Phi\left( u, h\right) \\
    \od{h}{t} &= \Psi\left( u, h\right) \,,
\end{aligned} \label{eqn:odesys}
\end{equation}
\end{linenomath}
where \( t \) is the time coordinate.
Let \( u^n \approx u(t^n)\) and \( h^n \approx h(t^n) \) be the numerical approximations to \( u \) and \( h \) at time \( t = t^n \).
Let \( \Delta t \) be a time-step such that \( t^{n+1} = t^n + \Delta t \).
Also, let \( t^{n+\nicefrac{1}{m}} = t^n + \frac{\Delta t}{m} \) for any positive integer \( m \).
Then, FB-RK(3,2) is given by
\begin{linenomath}
\begin{subequations}
\label{eqn:fbrk32}
\begin{align}
\begin{split}
    \bar{h}^{n+\nicefrac{1}{3}} &= h^n + \frac{\Delta t}{3} \Psi\left( u^n, h^n \right) \\
    \bar{u}^{n+\nicefrac{1}{3}} &= u^n + \frac{\Delta t}{3} \Phi\left( u^n, h^* \right) \\
    h^* &= \beta_1 \bar{h}^{n+\nicefrac{1}{3}} + (1-\beta_1) h^n
\end{split} \label{subeqn:fbrk32_stage1} \\
\nonumber \\
\begin{split}
    \bar{h}^{n+\nicefrac{1}{2}} &= h^n + \frac{\Delta t}{2} \Psi\left( \bar{u}^{n+\nicefrac{1}{3}}, \bar{h}^{n+\nicefrac{1}{3}} \right) \\
    \bar{u}^{n+\nicefrac{1}{2}} &= u^n + \frac{\Delta t}{2} \Phi\left( \bar{u}^{n+\nicefrac{1}{3}}, h^{**} \right) \\
    h^{**} &= \beta_2 \bar{h}^{n+\nicefrac{1}{2}} + (1-\beta_2) h^n
\end{split} \label{subeqn:fbrk32_stage2} \\
\nonumber \\
\begin{split}
    h^{n+1} &= h^n + \Delta t \Psi\left( \bar{u}^{n+\nicefrac{1}{2}}, \bar{h}^{n+\nicefrac{1}{2}} \right) \\
    u^{n+1} &= u^n + \Delta t \Phi\left( \bar{u}^{n+\nicefrac{1}{2}}, h^{***} \right) \\
    h^{***} &= \beta_3 h^{n+1} + (1-2\beta_3) \bar{h}^{n+\nicefrac{1}{2}} + \beta_3 h^n \,.
\end{split} \label{subeqn:fbrk32_stage3}
\end{align}
\end{subequations}
\end{linenomath}
The weights \( \beta_1 \), \( \beta_2 \), and \( \beta_3 \) are called the forward-backward (FB) weights. 
The best choice for the FB-weights \( \beta_i \) is a primary concern of this work.

The choice to calculate data at time \( t^{n+\nicefrac{1}{3}} \) and \( t^{n+\nicefrac{1}{2}} \) in the first and second stages comes from the original RK3 scheme; the treatment of the thickness variable \( h \) in \eqref{eqn:fbrk32} shows the formulation of RK3.
FB-RK(3,2) can be approximately reduced to RK3 be taking \( \beta_1 = 0 \), \( \beta_2 = \frac{2}{3} \), and \( \beta_3 = 0 \).
This produces an approximation to RK3 as, while the first and third stages are exact, the FB averaged \( h^{**} = \frac{2}{3} \bar{h}^{n+\nicefrac{1}{2}} + \frac{1}{3} h^n \) data in the second stage is an approximation to \( \bar{h}^{n+\nicefrac{1}{3}} \).
The FB averages are chosen in this way so that in each stage, data is symmetrically distributed in time as preliminary experiments showed that this provided the best CFL performance.
In the first two RK stages, the FB averaging is done on data from the beginning and the end of the current time interval, i.e. in the first stage, we compute values for \( h \) and \( \u \) at time \( t^{n+\nicefrac{1}{3}} \), so the FB averaging is done with data at time \( t^n \) and time \( t^{n+\nicefrac{1}{3}} \).
In the third RK stage, we use data from the beginning, middle, and end of the current time interval.

As stated above, the original RK3 scheme is \( \O\left( (\Delta t)^3 \right) \). 
For any choice of the FB-weights, FB-RK(3,2) reduces to at least \( \O\left( (\Delta t)^2 \right) \) while greatly increasing the maximum allowable time-step for most problems.
A generalized local truncation error analysis is outside the scope of this paper, so we show this holds for a particular problem of interest in Section \ref{subsec:lte_analysis}.


\subsection{Linearized Shallow Water Equations}
\label{subsec:linearized_swes}

A primary concern of this work is the application of FB-RK(3,2) to the nonlinear SWEs on a rotating sphere. Application of a von Neumann type analysis requires a linear set of PDEs, and we consider the two-dimensional system with linear Coriolis and linearized advection operators to ensure the time-step optimization is conducted on a system that is a close approximation to the fully nonlinear equations.
The linearized SWEs centered about a nonzero mean flow are given by
\begin{linenomath}
\begin{equation}
\begin{aligned}
    &\pd{\tilde{\u}}{t} + f\left(\tilde{\U} + \tilde{\u}\right)^\perp + \tilde{\zeta}\tilde{\U}^\perp = -\tilde{\U} \cdot \nabla \tilde{\u} - g \nabla \tilde{h} \\
    &\pd{\tilde{h}}{t} + H(\nabla \cdot \tilde{\u}) + \nabla \cdot \left(\tilde{h}\tilde{\U}\right) = 0 \,,
\end{aligned} \label{eqn:nondimless_linswe_nzmf}
\end{equation}
\end{linenomath}
where \( \tilde{\u} = \left( \tilde{u}(x, y, t), \tilde{v}(x, y, t) \right) \) is a perturbation of the horizontal fluid velocity,
\( x \) and \( y \) are the spatial coordinates,
\( t \) is the time coordinate,
\( f = 2\Omega\sin\phi \) is the Coriolis parameter,
\( \Omega\) is the angular velocity of the rotating sphere and 
\( \phi \) denotes latitude,
\( \tilde{\U} = \left( \tilde{U}, \tilde{V} \right) \) is the constant horizontal mean fluid velocity,
\( \tilde{\u}^\perp = \left( -\tilde{v}, \tilde{u} \right) \),
\( \tilde{\zeta} = \frac{\partial \tilde{v}}{\partial x} - \frac{\partial \tilde{u}}{\partial y} \) is the vertical component of vorticity \( \nabla \times \tilde{\u} \),
\( g \) is the gravitational constant,
\( \tilde{h} = \tilde{h}(x, y, t) \) is a perturbation of the layer thickness,
and \( H \) is the (constant in space) fluid layer thickness when at rest.
Note that while \( \tilde{\u} \) and \( \tilde{h} \) depend on \( x \) and \( y \) throughout this work, we will often suppress this dependence in the notation.
The full derivation of these equations is given in Section \ref{subsec:linearization}.

We can write these in terms of dimensionless quantities by writing \( c = \sqrt{gH} \), \( \u = \left( u, v \right) =\frac{\tilde{\u}}{c} \), \( \U = \frac{\tilde{\U}}{c} \), and \( \eta = \frac{\tilde{h}}{H} \).
Then, the dimensionless, linearized SWEs centered about a non-zero mean flow are given by
\begin{linenomath}
\begin{equation}
\begin{aligned}
    &\pd{\u}{t} + f\left( \U + \u \right)^\perp + c\zeta\U^\perp = -c\U \cdot \nabla\u - c \nabla \eta \\
    &\pd{\eta}{t} + c\left( \nabla \cdot \u \right) + c\left(\nabla \cdot (\eta\U)\right) = 0 \,.
\end{aligned} \label{eqn:linswe_nzmf}
\end{equation}
\end{linenomath}
Observe that when the mean flow \( \U \) is taken to be zero, \eqref{eqn:linswe_nzmf} becomes recognizable as the ``standard'' linearized SWEs:
\begin{linenomath}
\begin{equation}
\begin{aligned}
    &\pd{\u}{t} + f\u^\perp = -c \nabla \eta \\
    &\pd{\eta}{t} + c(\nabla \cdot \u) = 0 \,. 
\end{aligned} \label{eqn:swe}
\end{equation}
\end{linenomath}


\subsection{Von Neumann Formulation}
\label{subsec:von_neumann_formulation}

In the analysis that follows, we assume that the system \eqref{eqn:linswe_nzmf} is discretized on a rectangular, staggered Arakawa C-grid, where \( \eta \) is computed at cell centers, \( u \) is computed at the left and right edges of cells, and \( v \) is computed at the top and bottom edges of cells \citep{arakawa1977}.
For brevity, we skip some details of calculations relating to the spatial discretization in the main text and instead provide them in Section \ref{subsec:c-grid_discretization}.

With the goal of using a von Neumann-type analysis to obtain an optimal choice of FB-weights, we apply a Fourier transform in \( x \) and \( y \) to the linearized SWEs \eqref{eqn:linswe_nzmf}.
This is equivalent to seeking solutions of the form
\begin{linenomath}
\begin{align*}
    u(x, y, t) &= \hat{u}(k, \ell, t)e^{ikx+i\ell y} \\
    v(x, y, t) &= \hat{v}(k, \ell, t)e^{ikx+i\ell y} \\
    \eta(x, y, t) &= \hat{\eta}(k, \ell, t)e^{ikx+i\ell y} \,,
\end{align*}
\end{linenomath}
where \( k \) is the wavenumber with respect to \( x \), and \( \ell \) is the wavenumber with respect to \( y \).
In doing this, spatial derivatives are now all of the form \( \pd{}{*} \left( e^{ikx+i\ell y} \right) \).
Assuming a rectangular C-grid, using centered differences these spatial derivatives become
\begin{linenomath}
\begin{align}
    \pd{}{x} \left( e^{ikx+i\ell y} \right) &= \frac{i 2 \sin\left(k\frac{\Delta x}{2}\right)}{\Delta x}e^{ikx+i\ell y} \label{eqn:dx_e} \\
    \pd{}{y} \left( e^{ikx+i\ell y} \right) &= \frac{i 2 \sin\left(\ell\frac{\Delta y}{2}\right)}{\Delta y}e^{ikx+i\ell y} \,. \label{eqn:dy_e}
\end{align}
\end{linenomath}
To compute the Coriolis terms, we need the values of \( u \) at \( v \)-points and the values of \( v \) at \( u \)-points; call these \( u_v \) and \( v_u \) respectively.
On a rectangular C-grid, we achieve this with unweighted four-point averages of the nearest data, which gives
\begin{linenomath}
\begin{align}
    u_v(x, y, t) &= \cos\left(k\frac{\Delta x}{2}\right)\cos\left(\ell\frac{\Delta y}{2}\right)\hat{u}(k, \ell, t)e^{ikx+i\ell y} \label{eqn:coriolis_u} \\
    v_u(x, y, t) &= \cos\left(k\frac{\Delta x}{2}\right)\cos\left(\ell\frac{\Delta y}{2}\right)\hat{v}(k, \ell, t)e^{ikx+i\ell y} \,. \label{eqn:coriolis_v}
\end{align}
\end{linenomath}
The details of these calculations are given in Section \ref{subsec:c-grid_discretization}.
Set \( K = 2 \sin\left(k\frac{\Delta x}{2}\right) \), \( L = 2 \sin\left(\ell\frac{\Delta y}{2}\right) \), and \( \psi = f\cos\left(k\frac{\Delta x}{2}\right)\cos\left(\ell\frac{\Delta y}{2}\right) \). 
Then, accounting for the Fourier transform and the C-grid spatial discretization, \eqref{eqn:linswe_nzmf} becomes
\begin{linenomath}
\begin{equation}
\begin{aligned}
    \pd{\hat{u}}{t} &= fV + \psi\hat{v} - \left(iU\frac{cK}{\Delta x} + iV\frac{cL}{\Delta y}\right)\hat{u} - i\frac{cK}{\Delta x}\hat{\eta} \\
    \pd{\hat{v}}{t} &= -fU - \psi\hat{u} - \left(iU\frac{cK}{\Delta x} + iV\frac{cL}{\Delta y}\right)\hat{v} - i\frac{cL}{\Delta y}\hat{\eta} \\
    \pd{\hat{\eta}}{t} &= -\left( i\frac{cK}{\Delta x}\hat{u} + i\frac{cL}{\Delta y}\hat{v} \right) - \left( iU\frac{cK}{\Delta x} + iV\frac{cL}{\Delta y}\right)\hat{\eta} \,.
\end{aligned} \label{eqn:vn_linswe_nzmf}
\end{equation}
\end{linenomath}

Let \( \nu_x = c\frac{\Delta t}{\Delta x} \) and \( \nu_y = c\frac{\Delta t}{\Delta y} \) be the Courant numbers with respect to \( x \) and \( y \) respectively.
Also, set \( \phi = \Delta t\, \psi \).
Then, we apply FB-RK(3,2) to \eqref{eqn:vn_linswe_nzmf}:
\begin{linenomath}
\begin{subequations}
\label{eqn:vn_linswe_nzmf_fbrk32}
\begin{align}
\begin{split}
    \bar{\hat{\eta}}^{n+1/3} &= \hat{\eta}^n - \frac{1}{3}\left( iK\nu_x \hat{u}^n + iL\nu_y \hat{v}^n + \left( iUK\nu_x + iVL\nu_y \right)\hat{\eta}^n \right) \\
    \bar{\hat{u}}^{n+1/3} &= \hat{u}^n + \frac{1}{3}\left( \Delta t fV + \phi \hat{v}^n - \left( iUK\nu_x + iVL\nu_y \right)\hat{u}^n - iK\nu_x \hat{\eta}^* \right) \\
    \bar{\hat{v}}^{n+1/3} &= \hat{v}^n + \frac{1}{3}\left( -\Delta t fU - \phi \hat{u}^n - \left( iUK\nu_x + iVL\nu_y \right)\hat{v}^n - iK\nu_x \hat{\eta}^* \right) \\
    \hat{\eta}^* &= \beta_1\bar{\hat{\eta}}^{n+1/3} + (1-\beta_1)\hat{\eta}^n
\end{split} \label{subeqn:vn_linswe_nzmf_fbrk32_stage1} \\
\nonumber \\
\begin{split}
    \bar{\hat{\eta}}^{n+1/2} &= \hat{\eta}^n - \frac{1}{2}\left( iK\nu_x \bar{\hat{u}}^{n+1/3} + iL\nu_y \bar{\hat{v}}^{n+1/3} + \left( iUK\nu_x + iVL\nu_y \right)\bar{\hat{\eta}}^{n+1/3} \right) \\
    \bar{\hat{u}}^{n+1/2} &= \hat{u}^n + \frac{1}{2}\left( \Delta t fV + \phi \bar{\hat{v}}^{n+1/3} - \left( iUK\nu_x + iVL\nu_y \right)\bar{\hat{u}}^{n+1/3} - iK\nu_x \hat{\eta}^{**} \right) \\
    \bar{\hat{v}}^{n+1/2} &= \hat{v}^n + \frac{1}{2}\left( -\Delta t fU - \phi \bar{\hat{u}}^{n+1/3} - \left( iUK\nu_x + iVL\nu_y \right)\bar{\hat{v}}^{n+1/3} - iK\nu_x \hat{\eta}^{**} \right) \\
    \hat{\eta}^{**} &= \beta_2\bar{\hat{\eta}}^{n+1/2} + (1-\beta_2)\hat{\eta}^n
 \end{split} \label{subeqn:vn_linswe_nzmf_fbrk32_stage2} \\
 \nonumber \\
 \begin{split}
    \hat{\eta}^{n+1} &= \hat{\eta}^n - \left( iK\nu_x \bar{\hat{u}}^{n+1/2} + iL\nu_y \bar{\hat{v}}^{n+1/2} + \left( iUK\nu_x + iVL\nu_y \right)\bar{\hat{\eta}}^{n+1/2} \right) \\
    \hat{u}^{n+1} &= \hat{u}^n + \left( \Delta t fV + \phi \bar{\hat{v}}^{n+1/2} - \left( iUK\nu_x + iVL\nu_y \right)\bar{\hat{u}}^{n+1/2} - iK\nu_x \hat{\eta}^{***} \right) \\
    \hat{v}^{n+1} &= \hat{v}^n + \left( -\Delta t fU - \phi \bar{\hat{u}}^{n+1/2} - \left( iUK\nu_x + iVL\nu_y \right)\bar{\hat{v}}^{n+1/2} - iK\nu_x \hat{\eta}^{***} \right) \\
    \hat{\eta}^{***} &= \beta_3 \hat{\eta}^{n+1} + (1-2\beta_3)\bar{\hat{\eta}}^{n+1/2} + \beta_3\hat{\eta}^n \,.
\end{split} \label{subeqn:vn_linswe_nzmf_fbrk32_stage3}
\end{align}
\end{subequations}
\end{linenomath}

Note that \eqref{eqn:vn_linswe_nzmf_fbrk32} forms a linear system in \( \hat{u}^n \), \( \hat{v}^n \), and \( \hat{\eta}^n \), that is, there is a matrix \( G \) and column vector \( \mathbf{b} \) such that
\begin{linenomath}
\begin{equation}
    \hat{\w}^{n+1} = G \hat{\w}^n + \mathbf{b} \,, \label{eqn:amp_matrix}
\end{equation}
\end{linenomath}
where \( \hat{\w}^n = (\hat{u}^n, \hat{v}^n, \hat{\eta}^n)^T \).
It follows from \eqref{eqn:amp_matrix} that
\begin{linenomath}
\begin{equation}
    \hat{\w}^{n} = G^n \hat{\w}^0 + (I - G^n) (I - G)^{-1} \mathbf{b} \,, \label{eqn:amp_matrix_n}
\end{equation}
\end{linenomath}
where \( I \) is the identity matrix.
Note that in the above, the super-script on \( \hat{\w} \) is an index on the time-step, whereas the super-script on \( G \) is an exponent. 
The matrix \( G \) is called the amplification matrix; properties of von Neumann analysis tell us that a time-stepping method is stable if the eigenvalues of its amplification matrix reside within in the unit circle \citep{leveque2007}.
This property is fundamental to the analysis in Section \ref{subsec:cfl_optimization}.


\subsection{CFL Optimization}
\label{subsec:cfl_optimization}

With FB-RK(3,2) applied to the von Neumann formulation of the linearized SWEs \eqref{eqn:vn_linswe_nzmf_fbrk32} and the corresponding amplification matrix \( G \) from \eqref{eqn:amp_matrix}, we can more clearly state the goal of the analysis presented in this section.
We are interested in selecting the FB-weights \( \beta_1 \), \( \beta_2 \), and \( \beta_3 \) in FB-RK(3,2) so that the Courant numbers \( \nu_x = c\frac{\Delta t}{\Delta x} \) and \( \nu_y = c\frac{\Delta t}{\Delta y} \) can be taken as large as possible while retaining stability of the method.

The eigenvalues of the amplification matrix \( G \) depend on the Courant numbers \( \nu_x \) and \( \nu_y \) and on \( \beta_1 \), \( \beta_2 \), and \( \beta_3 \). 
To simplify the analysis, we will assume that \( \nu_x = \nu_y = \nu \) is the Courant number (this is equivalent to assuming that the spatial discretization is on a grid consisting of sqaure cells).
Then, we can write \( G = G(\nu, \beta_1, \beta_2, \beta_3) \).
For a given choice of \( \beta_1 \), \( \beta_2 \), and \( \beta_3 \), let \( \nu^\text{max}_{(\beta_1, \beta_2, \beta_3)} \) be the largest value of the Courant number for which FB-RK(3,2) is stable, i.e. largest value of the Courant number for which the eigenvalues of the amplification matrix \( G \) reside within the unit circle.
Symbolically, \( \nu^\text{max}_{(\beta_1, \beta_2, \beta_3)} \) is the largest number such that
\begin{linenomath}
\begin{equation}
    \abs{ \lambda_\text{max}\left( G(\nu, \beta_1, \beta_2, \beta_3) \right) } \leq 1 \qquad \text{for all } \nu \in \left(0, \nu^\text{max}_{(\beta_1, \beta_2, \beta_3)} \right] \,, \label{eqn:nu_max_defn}
\end{equation}
\end{linenomath}
where \( \lambda_\text{max}(G) \) is the largest eigenvalue, in absolute value, of the matrix \( G \).

While the matrix \( G \) is too complex to analyze directly, we can use numerical optimization tools to find optimal FB-weights.
That is, we can formulate this problem explicitly as an optimization problem with an appropriate cost function and constraints, then use numerical optimization algorithms to search the space of admittable FB-weights.
We can maximize \( \nu^\text{max} \) by minimizing an appropriate cost function
\begin{linenomath}
\begin{equation}
    C = C(\beta_1, \beta_2, \beta_3) \,,
\end{equation}
\end{linenomath}
subject to the constraint that the eigenvalues of \( G \) satisfy \eqref{eqn:nu_max_defn}.

The particular numerical optimization scheme we employ here is the simplicial homology global optimization (SHGO) algorithm developed by \citet{endres2018}, implemented in the Python programming language via the open-source Python library SciPy \citep{virtanen2020}.
SHGO is suitable for so-called black-box optimization problems wherein the cost function can be non-smooth and discontinuous.
Further, SHGO is guaranteed to converge to a global optimum as the number of iterations increases. 
As such, SHGO is an appropriate choice for our purposes, where the eigenvalues of the matrix \( G \) are sensitive to small changes in the FB-weights.

The matrix \( G \) also depends on additional parameters \( \Delta t\, f \), \( k \Delta x\), and \( \ell \Delta y\) as can be seen in \eqref{eqn:vn_linswe_nzmf_fbrk32}.
For the purposes of the numerical optimizations we perform here, we set \( \Delta t\, f = 10^{-2} \) since taking \( f = 10^{-4} \) represents the value of the Coriolis parameter at mid-latitudes and we expect a FB-RK(3,2) to be capable of a time-step on the order of \( 10^2 \) in the nonlinear SWE test cases described in Section \ref{sec:experiments}.
Preliminary tests showed that the choice of value for \( \Delta t\, f \) had little-to-no effect on the optimizations in practice.
We take \( k \Delta x = \ell \Delta y = \pi \), which corresponds to the case of grid-scale waves.
This choice is informed by the fact that instability often comes from grid-scale processes. 
That is, we perform the optimization on grid-scale waves, which represents a worst-case with the hope that this provides the best overall stability.

\subsubsection{The Choice of Cost Functions}
\label{subsubsec:choice_of_cost_function}

Given a particular choice of \( \beta_1 \), \( \beta_2 \), and \( \beta_3 \), we can easily estimate the value of \( \numaxbetas \) to arbitrary precision. Our goal is to maximize this parameter, so we seek to minimize the cost function
\begin{linenomath}
\begin{equation}
    C_1(\beta_1, \beta_2, \beta_3) = \frac{1}{\numaxbetas} \,. \label{eqn:cost_fun_c1}
\end{equation}
\end{linenomath}

In the special case where the mean flow \( \U = (U, V) \) is zero, we additionally consider a second choice of cost function. 
When \( \U = \mathbf{0} \) the column vector \( \mathbf{b} \) in \eqref{eqn:amp_matrix} is also zero, so \eqref{eqn:amp_matrix_n} simplifies to 
\begin{linenomath}
\begin{equation}
    \hat{\w}^{n} = G^n \hat{\w}^0 \,. \label{eqn:amp_matrix_n_zmf}
\end{equation}
\end{linenomath}
In this case, it is easy to show that the numerical amplification matrix \( G \)  approximates \( \tilde{G} \coloneqq e^{A\Delta t} \), where
\begin{linenomath}
\begin{equation}
    A = \begin{pmatrix}
        0 & f & -ick \\ -f & 0 & -ic\ell \\ -ick & -ic\ell & 0
    \end{pmatrix} \,. \label{eqn:a}
\end{equation}
\end{linenomath}
The matrix \( \tilde{G} \) is the analytical counterpart to the numerical amplification matrix \( G \) that arises from solving a spatially continuous version of \eqref{eqn:vn_linswe_nzmf}.
In short, if we assume that our time-stepping method, which is determined by the matrix \( G \), produces the exact solution at time \( t = \Delta t \), then \( G = \tilde{G} \).
We give the details of this in Appendix \ref{subsec:derivation_of_Gtilde}.

Since we wish for \( G \) to approximate \( \tilde{G} \) for a given choice of FB-weights, we aim to minimize a second choice of cost function:
\begin{linenomath}
\begin{equation}
    C_2(\beta_1, \beta_2, \beta_3) = \frac{1}{\numaxbetas} + \int_0^\frac{\pi}{6} \norm{\tilde{G}(\nu) - G(\nu, \beta_1, \beta_2, \beta_3)} \ \text{d}\nu \,, \label{eqn:cost_fun_c2}
\end{equation}
\end{linenomath}
where \( \norm{\cdot} \) is the Frobenius norm
\begin{linenomath}
\begin{equation}
        \norm{A} = \sqrt{\sum_{i,j}(a_{ij})^2} \,. \label{eqn:frob_norm}
\end{equation}
\end{linenomath}
This choice of cost function seeks to maximize \( \numax \) while also producing a scheme that has favorable dissipation and dispersion behaviour.
The choice of norm in the integral term of \( C_2 \) was chosen as it is quick to compute during the optimization process.
There is a certain freedom in the choice of limits of integration for the integral term in \( C_2 \); this term requires that the difference between \( G \) and \( \tilde{G} \) is minimal as for small values of \( \nu \), which correspond to well-resolved waves.
In essence, while \( C_1 \) only considers the CFL performance of the resulting method, \( C_2 \) does this while also considering that \( G \) should accurately approximate its analytical counterpart \( \tilde{G} \).

In summary, we consider two choices of cost function given by \( C_1 \) in \eqref{eqn:cost_fun_c1} and \( C_2 \) in \eqref{eqn:cost_fun_c2}.
The function \( C_1 \) is valid in any case, while \( C_2 \) applies only when the mean flow \( \U \) is taken to be zero.

\subsubsection{Optimization Results}
\label{subsubsec:optimization_results}

In the context of the dimensionless equations \eqref{eqn:linswe_nzmf}, the parameter \( \abs{\U} \) is called the Froude number.
We perform the optimization described above five total times, considering different values of the Froude number \( \abs{\U} \) that appear in oceanic and atmospheric flows.

\begin{table}
    \centering
    \renewcommand\arraystretch{1.5}
    \begin{tabular}{cc|cc}
         &  & \( \beta_1, \beta_2, \beta_3 \) & \( \numaxbetas \) \\
        \hline
        \multirow{2}*{\( \abs{\U} = \mathbf{0} \)} & \( C_1 \) & 0.500, 0.500, 0.344 & 1.767 \\
        & \( C_2 \) & 0.516, 0.532, 0.331 & 1.804 \\
        \( \abs{\U} = \mathbf{0.05} \) & \( C_1 \) & 0.531, 0.531, 0.313 & 1.319 \\
        \( \abs{\U} = \mathbf{0.15} \) & \( C_1 \) & 0.359, 0.578, 0.234 & 1.025 \\
        \( \abs{\U} = \mathbf{0.25} \) & \( C_1 \) & 0.656, 0.938, 0.188 & 0.853
    \end{tabular}
    \caption{
        Results of optimizing FB-RK(3,2) for the given cost functions, rounded to three decimal places. Note that the reported values of \( \numaxbetas \) are the maximal Courant numbers in the context of a problem consisting only of grid-scale waves, i.e. the problem for which FB-RK(3,2) was optimized.
    }
    \label{tbl:optmization_results}
\end{table}

As noted above, the SHGO algorithm used here is guaranteed to converge to a globally optimal solution if run for enough iterations.
As our problem is relatively complex and computing resources finite, it is more likely in practice that the method will convergence to a local extremum.
For example, consider the first two rows of Table \ref{tbl:optmization_results} and notice that the value of \( \numaxbetas \) given by optimizing \( C_2 \) is greater than that given by optimizing \( C_1 \).
If these were globally optimal solutions, we would have that \( \numax_{(0.500, 0.500, 0.344)} \geq \numax_{(0.516, 0.532, 0.331)} \) since \( C_1 \) consists only of a term that asks for a maximal \( \numaxbetas \), while \( C_2 \) is the sum of the same term with a nonnegative term that asks for \( G \) to approximate \( \tilde{G} \).
Despite the problem being sufficiently complex that globally optimal solutions are difficult to obtain, the results from Section \ref{sec:experiments} show that the  choices of FB-weights from Table \ref{tbl:optmization_results} are extremely efficient in practice.

A similar Runge-Kutta-based method using forward-backward averaging in the momentum equation is developed by \citet{shchepetkin2005} (from here on denoted ShMc).
In ShMc, the efficiency of their optimized FB method is reported in terms of the parameter \( \alpha \coloneqq ck\Delta t \), by considering the maximal value \( \alpha_\text{max} \) for which the method is stable.
Translating this into the notation used here, we obtain
\begin{linenomath}
\begin{align}
    \alpha &= ck\Delta t \nonumber \\
    &= \nu k\Delta x \,. \label{eqn:shchepetkin_alpha_tanslation}
\end{align}
\end{linenomath}
For their optimal RK2-based FB method, ShMc reports that \( (\nu k\Delta x)_\text{max} \approx 2.141 \); since their scheme uses two RK stages, this translates to what we refer to as an effective CFL number of \( 2.141 / 2 \approx 1.071  \).
Here, we have taken \( k\Delta x = \pi \) and have \( \numax = 1.804 \) in the best case, so \( (\nu k\Delta x)_\text{max} = 1.804\pi \approx 5.667 \).
FB-RK(3,2) uses three RK stages, so we obtain an effective CFL number of  \( 5.667 / 3 \approx 1.889 \), an increase of approximately 1.764 times over ShMc's RK2 FB scheme.
Note that ShMc's \( (\nu k\Delta x)_\text{max} \) comes from their FB scheme as applied to a one-dimensional wave system analogous to \eqref{eqn:vn_1d_linswe}, whereas ours is obtained from the two-dimensional linearized SWEs.
The comparison between these two values is not one-to-one, but do serve as a general comparison of the efficiency of both methods. 

\begin{figure}
    \def\width{0.3\textwidth}
    \centering
    \begin{subfigure}{\width}
        \includegraphics[width=\textwidth]{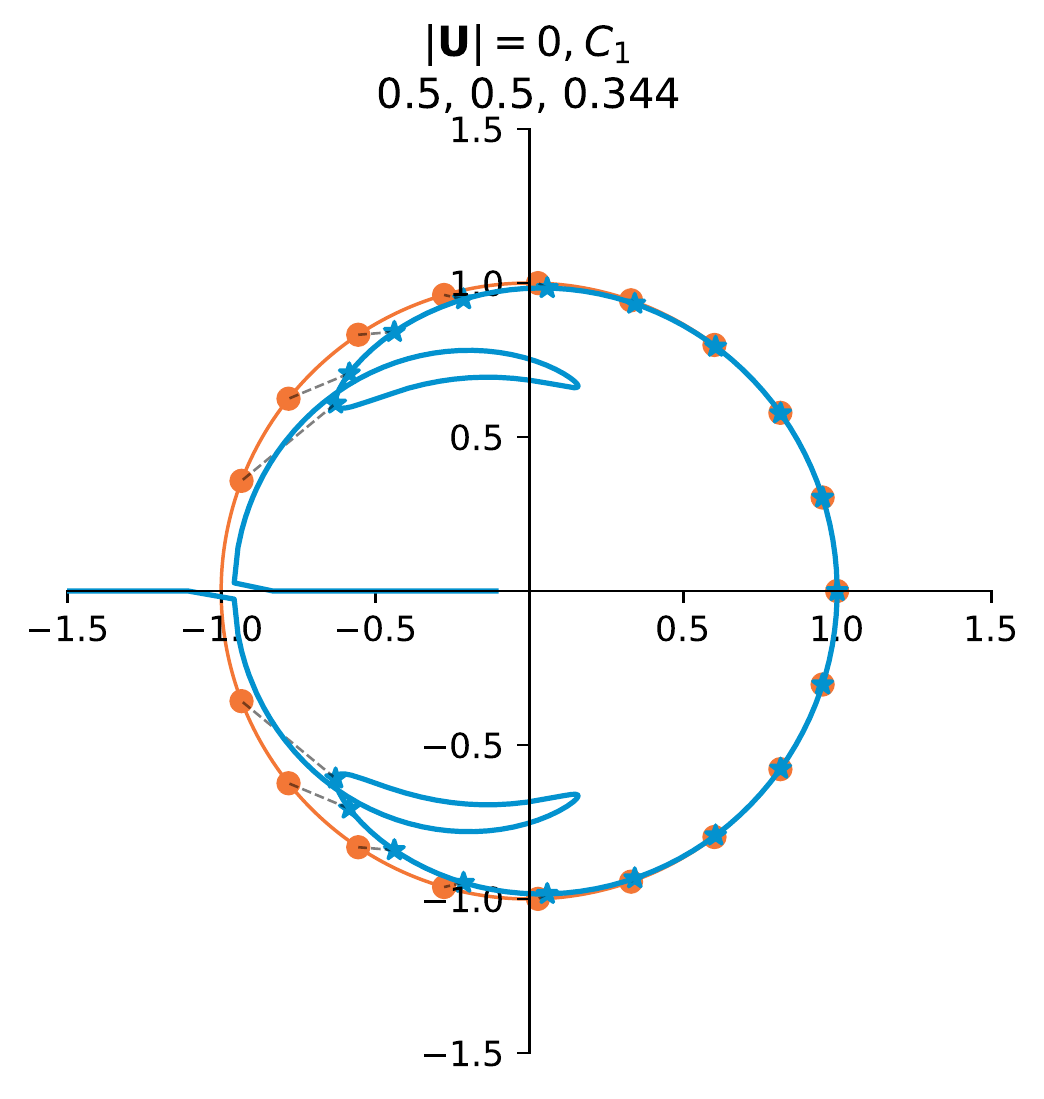}
        \caption{~}
        \label{subfig:U0C1}
    \end{subfigure}\begin{subfigure}{\width}
        \includegraphics[width=\textwidth]{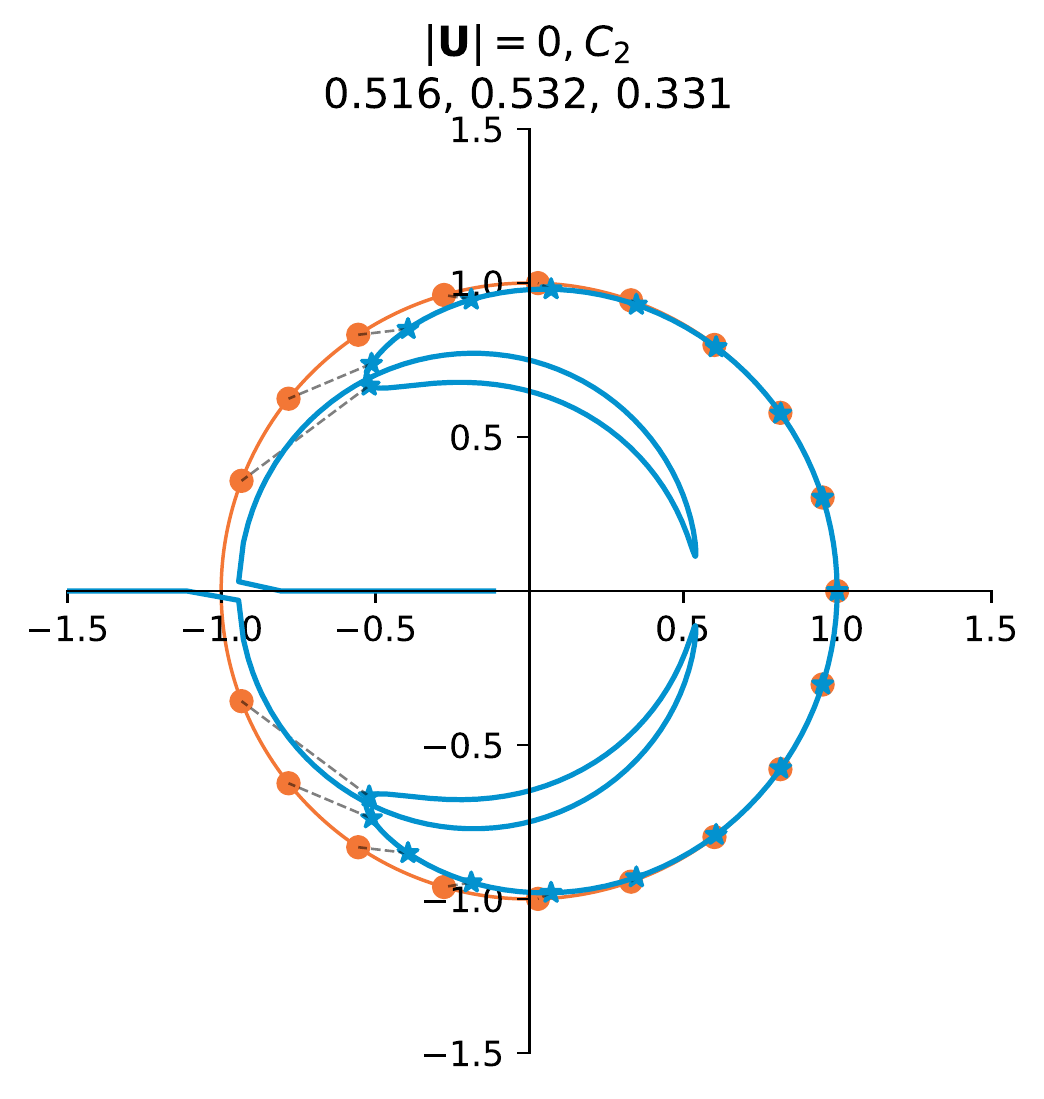}
        \caption{~}
        \label{subfig:U0C2}
    \end{subfigure}
    \begin{subfigure}{\width}
        \includegraphics[width=\textwidth]{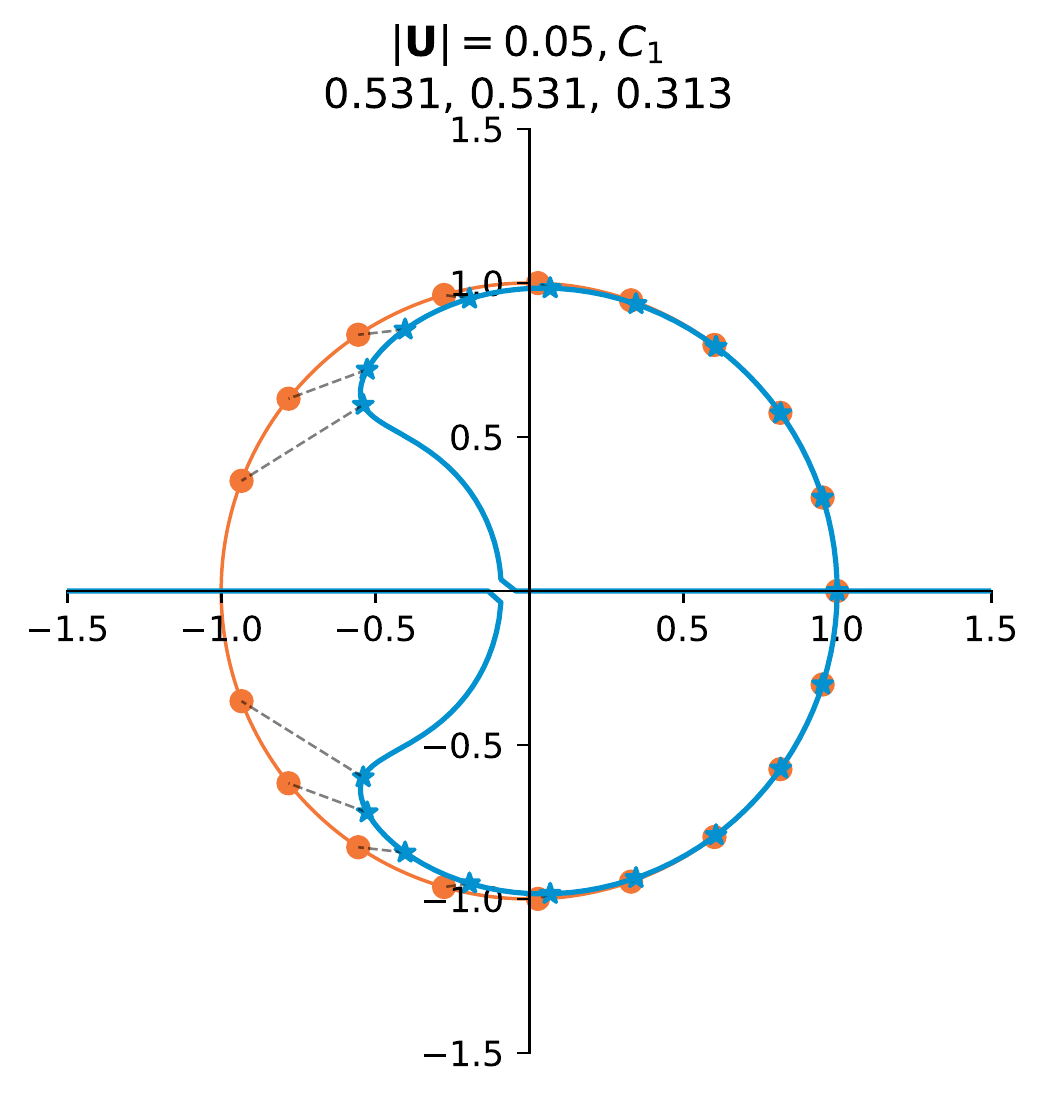}
        \caption{~}
        \label{subfig:U05C1}
    \end{subfigure}
    \begin{subfigure}{\width}
        \includegraphics[width=\textwidth]{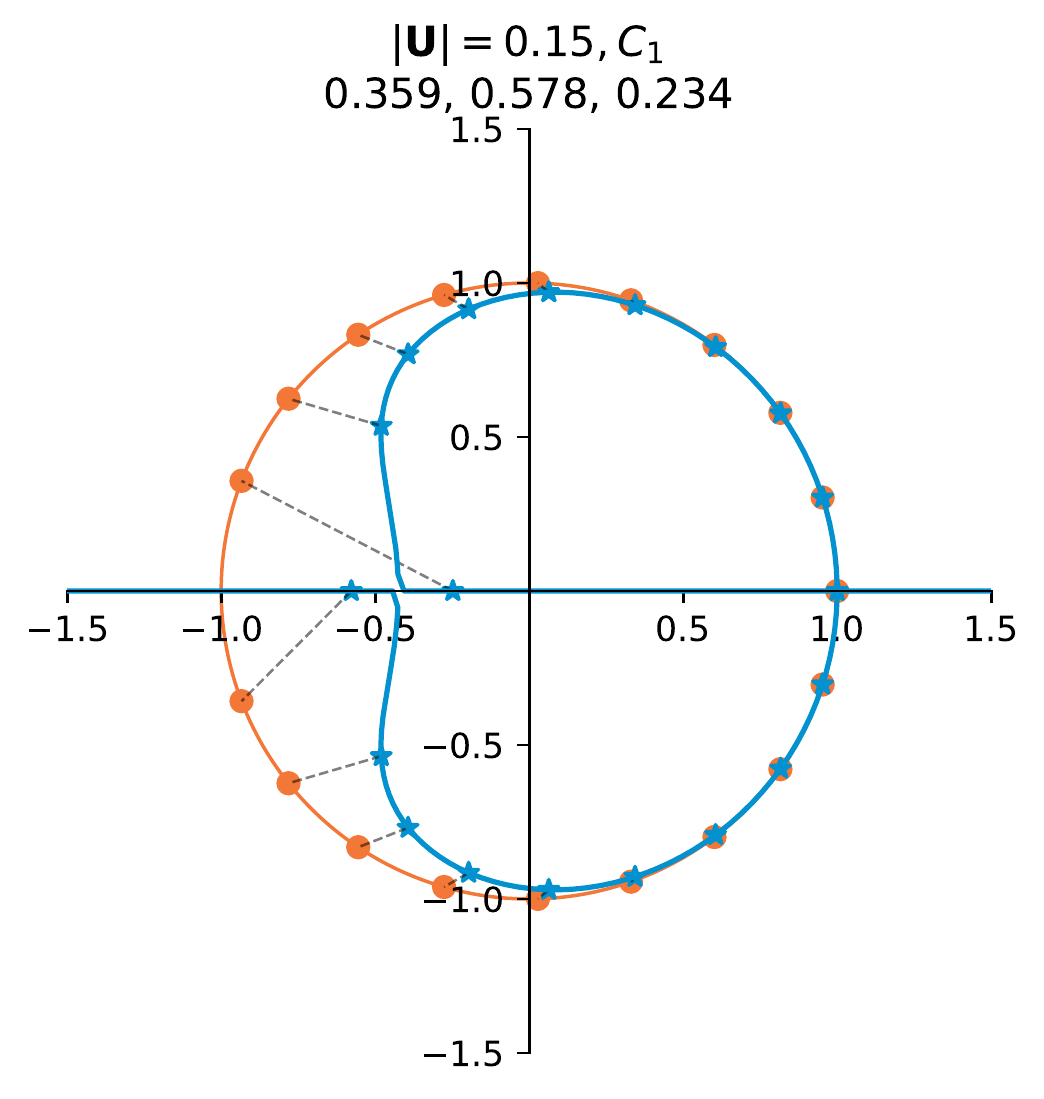}
        \caption{~}
        \label{subfig:U15C1}
    \end{subfigure}
    \begin{subfigure}{\width}
        \includegraphics[width=\textwidth]{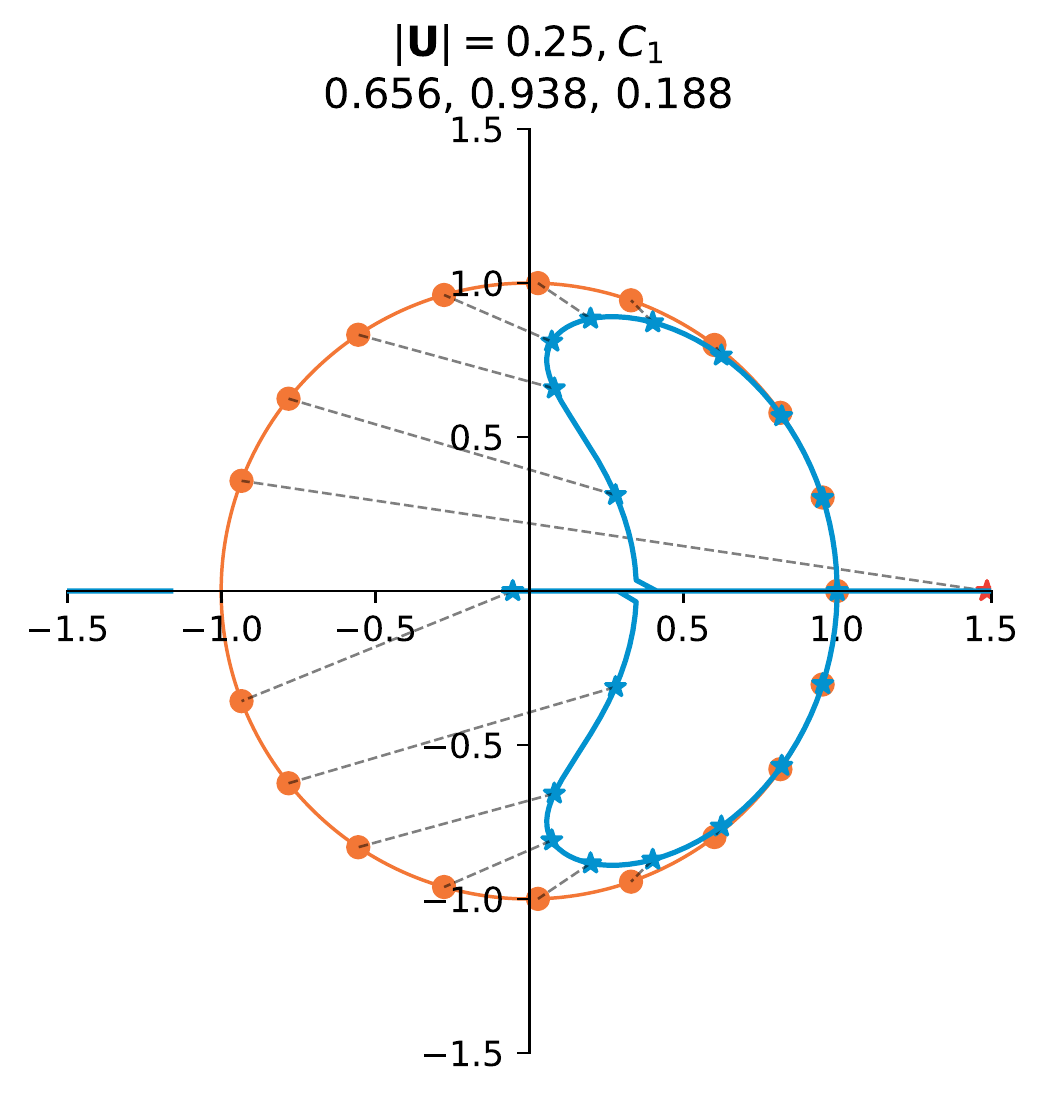}
        \caption{~}
        \label{subfig:U25C1}
    \end{subfigure}
    \caption{
        An illustration of dissipation and dispersion errors for FB-RK(3,2) with each choice of FB-weights from Table \ref{tbl:optmization_results}.
        Eigenvalues \( \lambda = e^{\pm i\tilde{K}\nu}\) of the matrix in the exact solution \eqref{eqn:vn_1d_linswe_general_soln} are given by orange dots.
        The angle these dots make with the positive real axis gives the corresponding value of \( \tilde{K}\nu \).
        Eigenvalues of the numerical amplification matrix \( G_0 \) as \( \tilde{K}\nu = k\Delta x\, \nu \) increases are given by the blue curve.
        Dashed lines connect exact eigenvalues to numerically determined counterparts, which gives a measure of dissipation and dispersion errors.
        Values of \( \tilde{K}\nu \) close to zero correspond to well resolved waves, while values approaching \( \pi \) correspond to (temporal) grid-scale waves.
    }
    \label{fig:eig_vals}
\end{figure}

To visualize FB-RK(3,2)'s treatment of gravity waves within the SWEs, we can apply the scheme to a simple wave system and analyze the resulting dissipation and dispersion errors (Figure \ref{fig:eig_vals}).
Applying a Fourier transform in space to the 1D linearized SWEs \eqref{eqn:1d_linswe} gives
\begin{linenomath}
\begin{equation}
\begin{aligned}
    \pd{\hat{\eta}}{t} &= -ick \hat{u} \\
    \pd{\hat{u}}{t} &= -ick \hat{\eta} \,.
\end{aligned} \label{eqn:vn_1d_linswe}
\end{equation}
\end{linenomath}
Choose a time-step \( \Delta t \) and let \( \tilde{K} = k\Delta x \) and \( \nu = c\frac{\Delta t}{\Delta x} \). 
The FB-RK(3,2) scheme for this system is then given by:
\begin{linenomath}
\begin{subequations}
\label{eqn:vn_1d_linswe_fbrk32}
\begin{align}
\begin{split}
    \bar{\eta}^{n+1/3} &= \eta^n - \frac{i \tilde{K}\nu}{3} u^n \\
    \bar{u}^{n+1/3} &= u^n - \frac{i \tilde{K}\nu}{3} \left( \beta_1\bar{\eta}^{n+1/3} + (1 - \beta_1)\eta^n \right)
\end{split} \label{subeqn:vn_1d_linswe_fbrk32_stage1} \\
\nonumber \\
\begin{split}
    \bar{\eta}^{n+1/2} &= \eta^n - \frac{i \tilde{K}\nu}{2} \bar{u}^{n+1/3} \\
    \bar{u}^{n+1/2} &= u^n - \frac{i \tilde{K}\nu}{2} \left( \beta_2 \bar{\eta}^{n+1/2} + (1 - \beta_2)\eta^n \right)
\end{split} \label{subeqn:vn_1d_linswe_fbrk32_stage2} \\
\nonumber \\
\begin{split}
    \eta^{n+1} &= \zeta^n - i \tilde{K}\nu \bar{u}^{n+1/2} \\
    u^{n+1} &= u^n -i \tilde{K}\nu \left( \beta_3\eta^{n+1} + (1-2\beta_3)\bar{\eta}^{n+1/2} + \beta_3\eta^n \right) \,.
\end{split} \label{subeqn:vn_1d_linswe_fbrk32_stage3}
\end{align}
\end{subequations}
\end{linenomath}
Equations \eqref{eqn:vn_1d_linswe_fbrk32} define a system of linear equations which in turn defines a numerical amplification matrix \( G_0 \) similar to that in \eqref{eqn:amp_matrix}.
One can show that the general solution to the system \eqref{eqn:vn_1d_linswe} is
\begin{linenomath}
\begin{equation}
    \begin{pmatrix} \eta(t) \\ u(t) \end{pmatrix} = \exp\left(t \begin{pmatrix} 0 & -ick \\ -ick & 0 \end{pmatrix}\right) \begin{pmatrix} \eta(0) \\ u(0) \end{pmatrix}, \label{eqn:vn_1d_linswe_general_soln}
\end{equation}
\end{linenomath}
and that eigenvalues of the matrix exponential above are \( \lambda = e^{\pm ick t} \). 
At time \( t = \Delta t \), these eigenvalues are \( \lambda = e^{\pm i\tilde{K}\nu} \).

Note that although we have used the Courant number \( \nu = c\frac{\Delta t}{\Delta x} \) in this discussion in an effort to use notation analogous to the rest of this section, the system considered here is continuous in space; we have that \( \tilde{K}\nu = k\Delta x \cdot c\frac{\Delta t}{\Delta x} = ck \Delta t \).
Values of \( \tilde{K}\nu \) close to zero correspond to well-resolved waves, while values approaching \( \pi \) correspond to (temporal) grid-scale waves.
Dissipation error can be read as the distance a numerical eigenvalue is from the unit circle, resulting in the corresponding wave being damped out of the system.
Dispersion error can be read as the difference in the angles made with the positive real axis between a numerically determined eigenvalue with its analytical counterpart, resulting in an error in the phase speed of the corresponding wave.
For all choices of FB-weights in Table \ref{tbl:optmization_results}, dissipation and dispersion errors for well resolved waves are very low (Figure \ref{fig:eig_vals}).
As waves become less well resolved, each of the five schemes exhibits different behavior.
In particular, note the behavior of the eigenvalues in Figure \ref{subfig:U05C1}; well resolved waves are handled well as expected, and as waves move towards grid-scale, the dissipation errors become very high, i.e. the blue curve becomes closer to the origin.
This means that waves that are not well-resolved are quickly damped out of the system, which may explain the strong performance of FB-RK(3,2) with the corresponding choices of FB-weights.

Note that the paradigm in which the plots of Figure \ref{fig:eig_vals} exist is different than the paradigm for which the FB-weights were optimized.
That is, the matrix \( G_0 \) comes from applying FB-RK(3,2) to a one-dimensional gravity wave system \eqref{eqn:vn_1d_linswe}, whereas the FB-weights presented in Table \ref{tbl:optmization_results} are optimized for the linearized SWEs \eqref{eqn:vn_linswe_nzmf}.
We present the plots in Figure \ref{fig:eig_vals} to show that the schemes optimized for stability in the linearized SWEs treat well-resolved gravity waves correctly.


\section{Numerical Experiments}
\label{sec:experiments}

Here, we present a number of numerical experiments that demonstrate the efficiency and accuracy of the CFL-optimized FB-RK(3,2) schemes presented in Section \ref{sec:optimization}.
These experiments compare SSPRK3 and RK3 to FB-RK(3,2) in terms of time-stepping efficiency, temporal convergence, and solution accuracy.
We consider a total of five different nonlinear shallow water test cases described below.
The fully nonlinear SWEs are given by \eqref{eqn:nonlinear_swe} in Section \ref{subsec:linearization}.

All test cases take place on a rotating aquaplanet with angular velocity \( \Omega = 7.292\times10^{-5} \) \( \text{s}^{-1} \) and a radius of 6371.22 km.
The sea-floor topography is uniform unless otherwise specified.
Each test case is implemented using a TRiSK spatial discretization \citep{ringler2010}, which is a finite volume-type spatial discretization made for unstructured, variable-resolution polygonal grids.
Balanced initial conditions for the BUJ and WTC test cases are generated following the procedure described in Section \ref{subsec:ICs}.
The code used to run all the experiments described in this section is open source and available on GitHub and Zenodo; see the data availability statement at the end of this document.

\begin{itemize}
    \item[] \textbf{Quasi-Linear Gravity Wave (QLW)}: A simple gravity wave traveling from the north pole to the south, then returning to the north, a process which takes approximately seven days of simulated time.
    The momentum advection terms in the momentum equation of \eqref{eqn:nonlinear_swe} are turned off, making the momentum equation purely linear.
    The thickness equation \eqref{eqn:nonlinear_swe} is formally nonlinear, but because the fluid depth is sufficiently large compared to any perturbation in the thickness, this behaves as though it is linear. 
    Therefore, we refer to this test case as quasi-linear.
    The fluid velocity is initialized to zero, while the thickness is initialized as a Gaussian bell curve centered at the north pole of the form \( e^{-100(x-\pi)^2 - 100y^2} + 500 \).

    \smallskip
    \item[] \textbf{Barotropically Unstable Jet (BUJ):} Originally formulated by \citet{galewsky2004}, this case consists of a strong zonal flow at mid-latitude with a small perturbation that causes the barotropic jet to become unstable. 
    This flow is driven by strongly nonlinear momentum advection tendencies -- physical processes that are of particular relevance to weather prediction and climate modeling.

    \smallskip
    \item[] \textbf{Williamson Test Case 2 (WTC2)}: A steady state solution to the nonlinear SWEs, one of the standard SWE test cases developed by \citet{williamson1992}.
    The case consists of a geostrophically balanced zonal flow that is constant for all time.

    \smallskip
    \item[] \textbf{Williamson Test Case 4 (WTC4)}: A forced nonlinear flow consisting of a translating low pressure center superimposed on a jet stream, one of the standard \citet{williamson1992} SWE test cases.

    \smallskip
    \item[] \textbf{Williamson Test Case 5 (WTC5)}: A zonal flow over an isolated mountain, one of the standard \citet{williamson1992} SWE test cases.
    Here, the initial flow is exactly as in WTC2, but the sea-floor topography is nonuniform, consisting of an isolated mountain at mid-latitude. This perturbation alters the geostrophic balance of the flow, triggering the development of advection-driven instabilities.
\end{itemize}


\subsection{CFL Performance}
\label{subsec:cfl_performance}

On all of the test cases described above, the optimized FB-RK(3,2) schemes from Section \ref{subsec:cfl_optimization} admit time-steps between 1.59 and 2.81 times as large as SSPRK3 or RK3.
This demonstrates the ability for the FB averaging employed by FB-RK(3,2) and the subsequent optimization of the FB weights to greatly increase the stability range of an existing method across multiple problems.
SSPRK3 is used as a baseline comparison here since it admits the smallest time-steps for each case out of the schemes considered here, and it also consists of three Runge-Kutta stages, costing the same number of floating-point operations as RK3 and FB-RK(3,2).

\begin{table}
    \centering
    \resizebox{\columnwidth}{!}{%
    \begin{tabular}{cc|ccccccc}
         &  & \multirow{4}*{SSPRK3} & \multirow{4}*{RK3} & \( \abs{\U} = 0,\ C_1 \) & \( \abs{\U} = 0,\ C_2 \) & \( \abs{\U} = 0.05,\ C_1 \) & \( \abs{\U} = 0.15,\ C_1 \) & \( \abs{\U} = 0.25,\ C_1 \) \\
         &  &  &  & 0.500 & 0.516 & 0.531 & 0.359 & 0.656 \\
         &  &  &  & 0.500 & 0.532 & 0.531 & 0.578 & 0.938 \\
        \multicolumn{2}{c|}{Test Case} &  &  & 0.344 & 0.331 & 0.313 & 0.234 & 0.188 \\
        \hline\hline
        \textbf{QLW} & 7 days &  &  &  &  &  &  &  \\
        \hline
        \multicolumn{2}{c|}{Max \( \Delta t \) (sec)} & 515 & 515 & 1445 & 1435 & 1115 & 960 & 820 \\
        \multicolumn{2}{c|}{ \( \nicefrac{\Delta t}{ \Delta t_\text{SSPRK3} } \) } &  & 1.00 & 2.81 & 2.79 & 2.17 & 1.86 & 1.59 \\
        \hline
        \textbf{BUJ} & 6 days &  &  &  &  &  &  &  \\
        \hline
        \multicolumn{2}{c|}{Max \( \Delta t \) (sec)} & 110 & 115 & 195 & 195 & 195 & 190 & 185 \\
        \multicolumn{2}{c|}{ \( \nicefrac{\Delta t}{ \Delta t_\text{SSPRK3} } \) } &  & 1.05 & 1.77 & 1.77 & 1.77 & 1.73 & 1.68 \\
        \hline
        \textbf{WTC2} & 5 days &  &  &  &  &  &  &  \\
        \hline
        \multicolumn{2}{c|}{Max \( \Delta t \) (sec)} & 220 & 230 & 355 & 360 & 370 & 410 & 360 \\
        \multicolumn{2}{c|}{ \( \nicefrac{\Delta t}{ \Delta t_\text{SSPRK3} } \) } &  & 1.05 & 1.61 & 1.64 & 1.68 & 1.86 & 1.64 \\
        \hline
        \textbf{WTC4} & 5 days &  &  &  &  &  &  &  \\
        \hline
        \multicolumn{2}{c|}{Max \( \Delta t \) (sec)} & 110 & 115 & 230 & 235 & 240 & 210 & 175 \\
        \multicolumn{2}{c|}{ \( \nicefrac{\Delta t}{ \Delta t_\text{SSPRK3} } \) } &  & 1.05 & 2.09 & 2.14 & 2.18 & 1.91 & 1.59 \\
        \hline
        \textbf{WTC5} & 15 days &  &  &  &  &  &  &  \\
        \hline
        \multicolumn{2}{c|}{Max \( \Delta t \) (sec)} & 145 & 145 & 270 & 275 & 310 & 275 & 235 \\
        \multicolumn{2}{c|}{ \( \nicefrac{\Delta t}{ \Delta t_\text{SSPRK3} } \) } &  & 1.00 & 1.86 & 1.9 & 2.14 & 1.90 & 1.62 \\
    \end{tabular}}
    \caption{
        CFL performance of FB-RK(3,2) with optimal FB-weights from Table \ref{tbl:optmization_results} compared to SSPRK3 and RK3 on various test cases.
        The recorded values for \( \Delta t \) represent the largest time-step for which the given test case could run without becoming unstable for the given duration.
        For practical reasons, time-steps were found as multiples of five.
    }
    \label{tab:cfl_performance}
\end{table}

The largest increases in the admittable time-step are in the QLW test case, with the largest increase at 2.81 times in the best case (Table \ref{tab:cfl_performance}).
For practical reasons, we used a nonlinear SWE model to run the above test cases, but QWL was intentionally designed so that the model equations, while still formally nonlinear, approximate the behavior of the linear SWEs.
That is, QLW provides a test case as close as possible to the paradigm for which FB-RK(3,2) was optimized.
As noted, the comparison is not one-to-one.
FB-RK(3,2) was optimized for the linearized SWEs on a square C-grid, while QLW uses quasi-linear SWEs and a hexagonal C-grid.
Nonetheless, the improvement in CFL performance is remarkable and clearly demonstrates the efficacy of the optimization performed in Section \ref{subsec:cfl_optimization}.

The remaining four test cases are fully nonlinear, covering a range of standard shallow water test cases.
In these, FB-RK(3,2) admits time-steps between 1.59 and 2.18 times larger than SSPRK3.
As SSPRK3 consists of the same number of Runge-Kutta stages as FB-RK(3,2), this directly translates to an increase of overall efficiency of between 1.59 and 2.18.
The increases in admittable time-step in these nonlinear cases are slightly less than those in the QLW case simply because these problems include fully nonlinear motions that cannot be accounted for in the optimization process used here; because our optimization relies on properties of von Nuemann stability analysis, we are limited to optimizing on a linear set of equations.
However, Table \ref{tbl:optmization_results} seems to suggest that optimizing for the linearized SWEs provides substantial benefits even for the nonlinear SWEs.

Temporarily ignoring the QLW case as it is not particularly relevant for real-world applications, we see that FB-RK(3,2) with the FB-weights obtained by optimizing \( C_1 \) with \( \abs{\U} = 0.05 \) admits the largest time step in each case except WTC2 where it instead admits the second largest time-step.
This points to \( ( \beta_1, \beta_2, \beta_3 ) = ( 0.531, 0.531, 0.313 ) \) as a robust option for a relatively broad span of problems.
However, it is almost certain though that the optimal choice of coefficient will be problem dependent.


\subsection{Temporal Convergence}
\label{subsec:temporal_convergence}

In Section \ref{subsec:lte_analysis}, we perform a simplified local truncation error analysis that shows that FB-RK(3,2) is second-order in time for most choices of FB-weights on the one dimensional linearized SWEs.
Here, we support the assertion that FB-RK(3,2) is second-order in time with numerical experiments.

\begin{figure}
    \centering
    \begin{subfigure}{0.75\textwidth}
        \includegraphics[width=\textwidth]{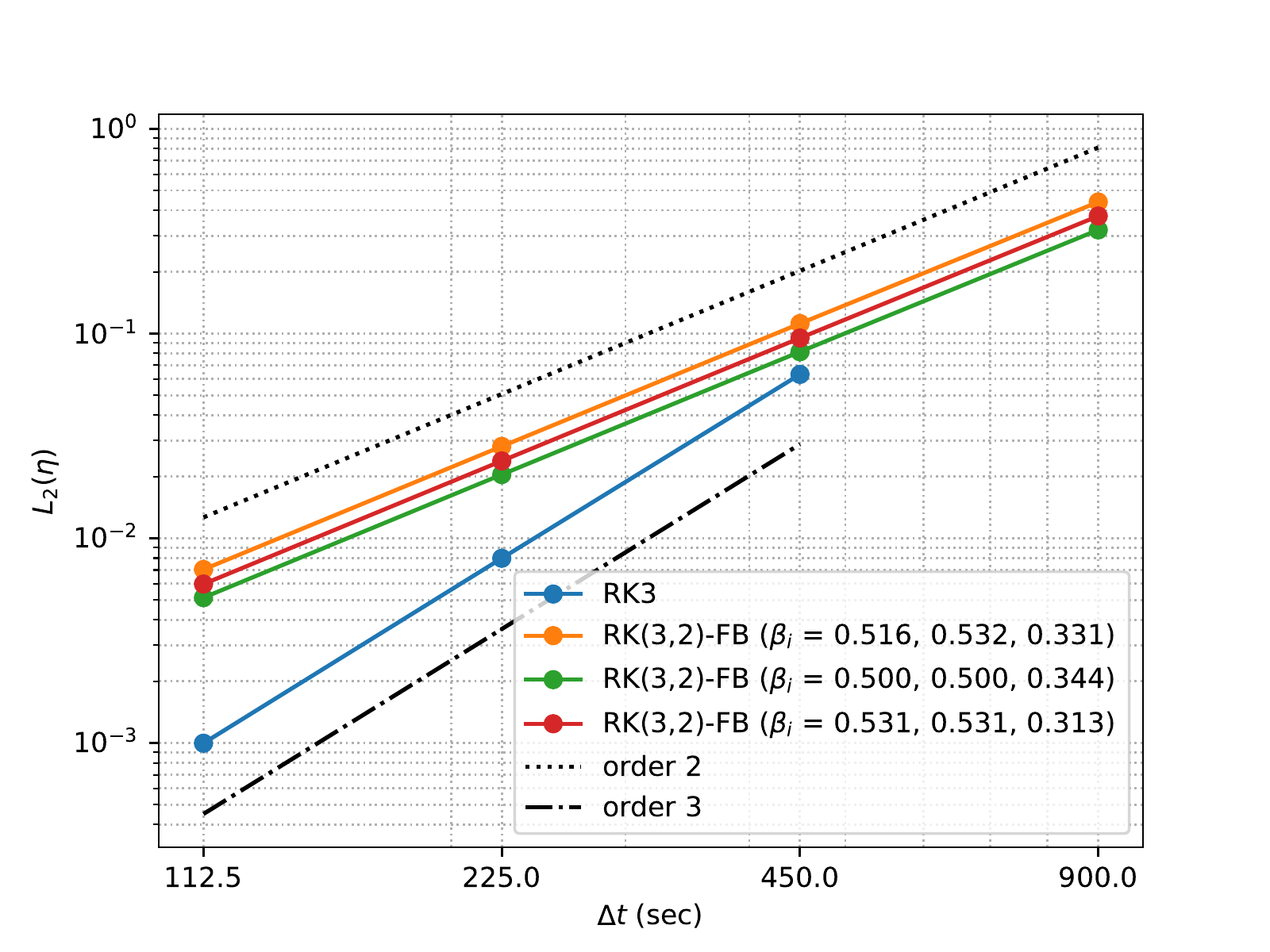}
        \caption{~}
        \label{subfig:temproal_convergence_eta}
    \end{subfigure} \\
    \begin{subfigure}{0.75\textwidth}
        \includegraphics[width=\textwidth]{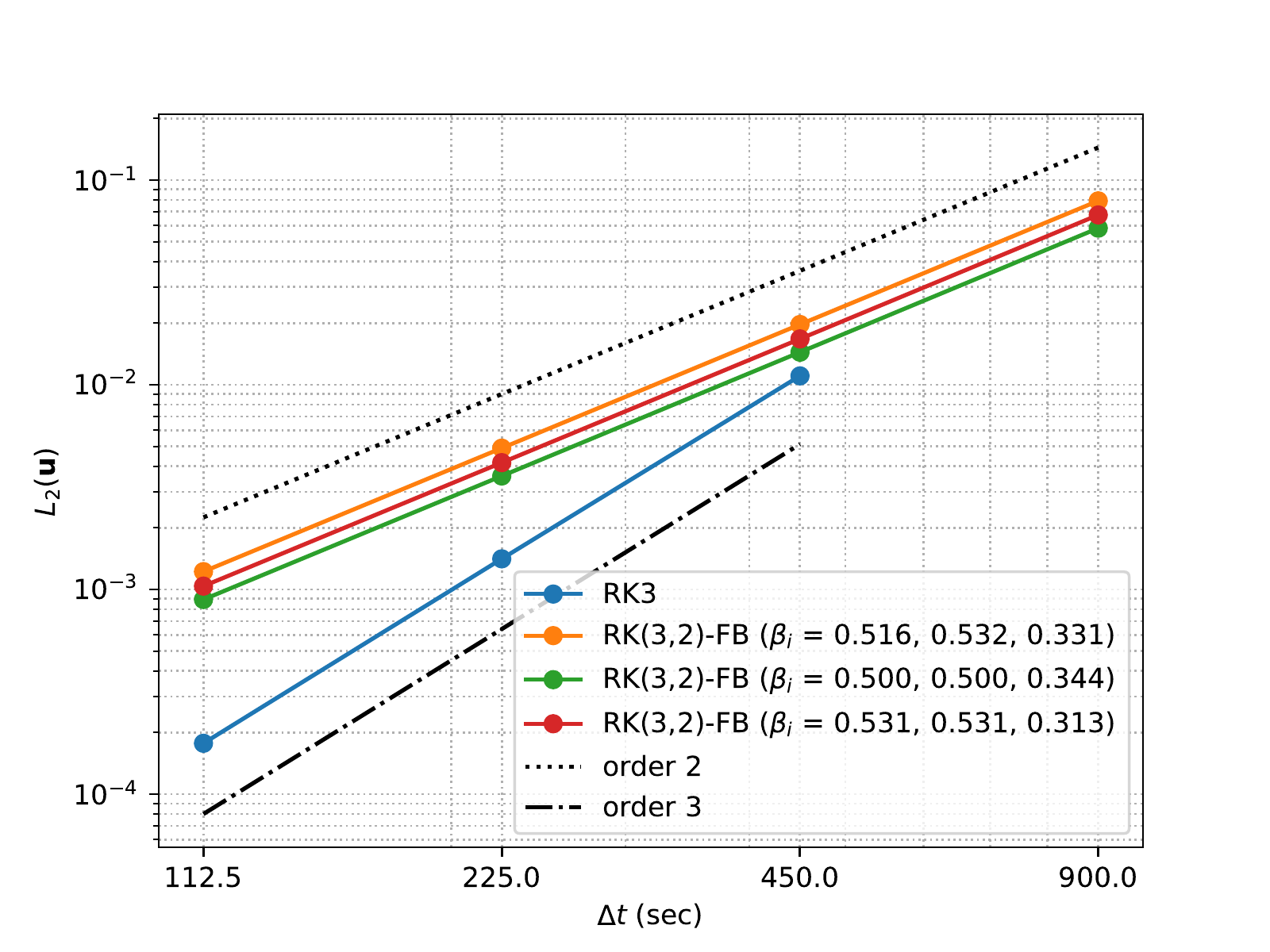}
        \caption{~}
        \label{subfig:temproal_convergence_u}
    \end{subfigure}
    \caption{
        Temporal convergence for FB-RK(3,2) for a selection of FB-weights.
        Convergence tests were performed using the QLW case, run for seven simulated days.
        \( L_2 \) errors are computed against a reference solution generated by the classical four-stage, fourth-order Runge-Kutta method (RK4) with a time-step of 10 s.
    }
    \label{fig:temproal_convergence}
\end{figure}

Each tested variation of FB-RK(3,2) is second-order in time (Figure \ref{fig:temproal_convergence}).
Differences in accuracy between the variations of FB-RK(3,2) are negligible, showing that the choice of optimized FB-weights has little effect on accuracy in this case. 
RK3, the scheme on which FB-RK(3,2) is based, is third-order as expected.
In effect, the addition of FB averaging to RK3 sacrifices a single order of accuracy for a greatly increased maximal Courant number.
In many applications, particularly in applications that incur lower-order errors from a spatial discretization or where right-hand-side terms are very expensive to evaluate, this trade-off is more than reasonable.


\subsection{Solution Quality}
\label{subsec:solution_quality}

We evaluate the solution quality of an optimal FB-RK(3,2) to that of SSPRK3 in two test cases of interest, BUJ and WTC5, comparing the relative vorticity produced by FB-RK(3,2) with \( ( \beta_1, \beta_2, \beta_3 ) = ( 0.531, 0.531, 0.313 ) \) to that produced by SSPRK3.
The strong stability preserving property of SSPRK3 makes it a good choice for a point of comparison; if FB-RK(3,2) can produce a solution of comparable quality, we can be confident that it is a good choice for the SWEs and related problems.

\begin{figure}
    \def\width{0.75\textwidth}
    \centering
    \begin{subfigure}{\width}
        \includegraphics[width=\textwidth]{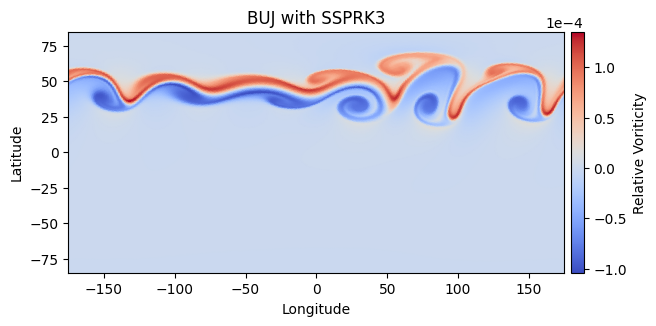}
        \caption{~}
        \label{subfig:buj_ssprk3}
    \end{subfigure}
    \begin{subfigure}{\width}
        \includegraphics[width=\textwidth]{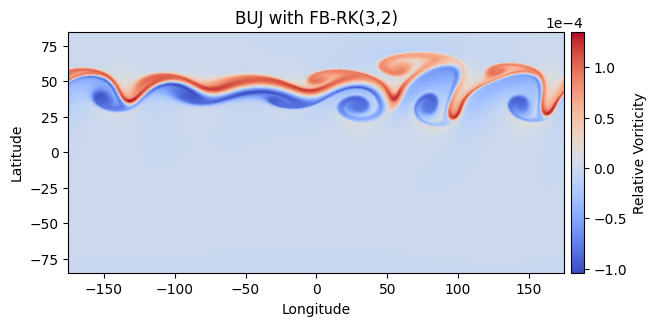}
        \caption{~}
        \label{subfig:buj_rk32fb}
    \end{subfigure}
    \begin{subfigure}{\width}
        \includegraphics[width=\textwidth]{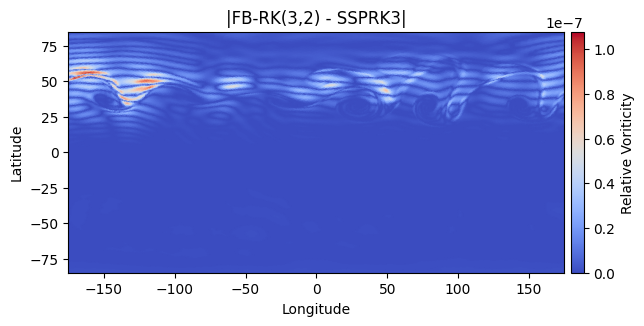}
        \caption{~}
        \label{subfig:buj_diff}
    \end{subfigure}
    \caption{
        Relative vorticity produced by SSPRK3 (Figure \ref{subfig:buj_ssprk3}) and FB-RK(3,2) with \( ( \beta_1, \beta_2, \beta_3 ) = ( 0.531, 0.531, 0.313 ) \) (Figure \ref{subfig:buj_rk32fb}) on the BUJ test case after six days of simulated time.
        SSPRK3 uses a time-step of 108 s and FB-RK(3,2) uses a time-step of 192 s.
        Figure \ref{subfig:buj_diff} shows the absolute difference between the two solutions.
        The plots themselves are lat/lon projections of the sphere.
    }
    \label{fig:buj}
\end{figure}

\begin{figure}
    \def\width{0.75\textwidth}
    \centering
    \begin{subfigure}{\width}
        \includegraphics[width=\textwidth]{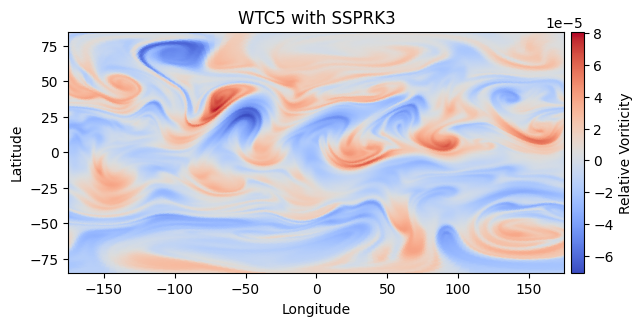}
        \caption{~}
        \label{subfig:wtc5_ssprk3}
    \end{subfigure}
    \begin{subfigure}{\width}
        \includegraphics[width=\textwidth]{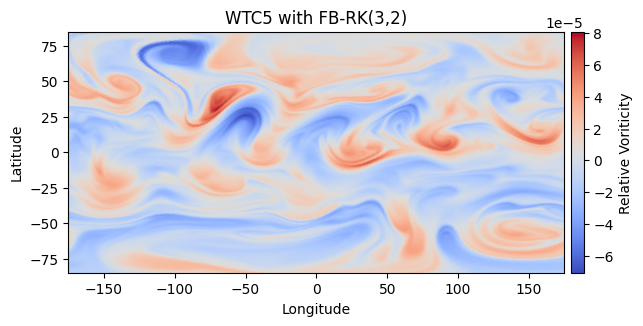}
        \caption{~}
        \label{subfig:wtc5_rk32fb}
    \end{subfigure}
    \begin{subfigure}{\width}
        \includegraphics[width=\textwidth]{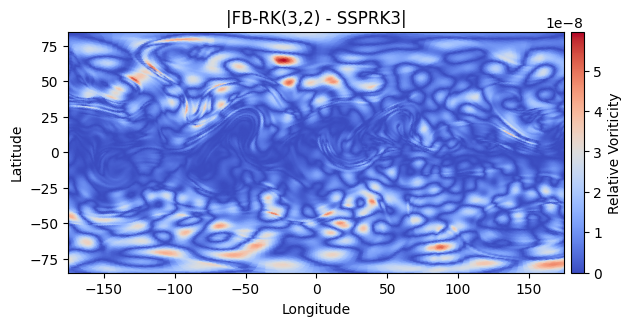}
        \caption{~}
        \label{subfig:wtc5_diff}
    \end{subfigure}
    \caption{
        Relative vorticity produced by SSPRK3 (Figure \ref{subfig:wtc5_ssprk3}) and FB-RK(3,2) with \( ( \beta_1, \beta_2, \beta_3 ) = ( 0.531, 0.531, 0.313 ) \) (Figure \ref{subfig:wtc5_rk32fb}) on the WTC5 test case after 50 days of simulated time.
        SSPRK3 uses a time-step of 135 s and FB-RK(3,2) uses a time-step of 288 s.
        Figure \ref{subfig:wtc5_diff} shows the absolute difference between the two solutions.
        The plots themselves are lat/lon projections of the sphere.
    }
    \label{fig:wtc5}
\end{figure}

In both cases, shown in Figure \ref{fig:buj} and Figure \ref{fig:wtc5}, FB-RK(3,2) produces solutions of comparable quality to SSPRK3 even with a significantly increased time-step.
In BUJ, the two schemes produce relative velocities that differ by at most \( 10^{-7} \).
In WTC5, this difference is on the order of \( 10^{-8} \).
We also ran both test cases using RK3 which produced plots, which we omit, that are visually identical to those given here.

FB-RK(3,2) uses a time-step 1.78 and 2.13 times larger that that used by SSPRK3 on BUJ and WTC5 respectively.
Despite taking about half the wall-clock time and computational resources, FB-RK(3,2) produces vorticities that are qualitatively equivalent to SSPRK3.
As noted in Section \ref{subsec:cfl_performance}, though FB-RK(3,2) is lower-order in time than SSPRK3, the produced solutions are qualitatively equivalent; likely partially due to lower-order spatial errors which dominate.
In operational configurations of climate-scale models, such spatial errors are generally dominant, so one can use a more economical time-stepping scheme like FB-RK(3,2) without sacrificing model quality.


\section{Conclusion}
\label{sec:conclusion}

With the goal of producing computationally efficient time-stepping schemes for the SWEs, we have presented a class of three-stage second-order Runge-Kutta schemes that use forward-backward averaging to advance the momentum equation.
The FB-weights in these FB-RK(3,2) schemes were then optimized using a von Neumann-type analysis on the shallow water equations linearized about a nonzero mean flow.
Even though it is impossible to perform a similar analysis and optimization on the fully nonlinear SWEs, the results of our linearized analysis clearly result in significant increases in efficiency when solving nonlinear test problems.

Optimized FB-RK(3,2) schemes are up to 2.8 times more computationally efficient than the RK3 scheme from which they were derived and a three-stage third-order strong stability preserving Runge-Kutta scheme (Table \ref{tab:cfl_performance}).
Despite being significantly more computationally economical that other schemes, FB-RK(3,2) produces solutions that are of comparable quality to other schemes while using a time-step approximately twice as long (Figure \ref{fig:buj} and \ref{fig:wtc5}).

Moving forward, we are interested in adapting FB-RK(3,2) to a local time-stepping (LTS) scheme similar to that used in \citep{dawson2013} or \citep{hoang2019}, thereby combining the benefits of the CFL performance of FB-RK(3,2) with the obvious performance implications of LTS schemes.
Another potentially valuable application of FB-RK(3,2) would be its use as a barotropic solver within a split-explicit, layered model wherein barotropic and baroclinic motions are separated and solved by different methods.
The barotropic equations are essentially the SWEs and our work here clearly demonstrates that FB-RK(3,2) is a good choice for solving such a system.



\clearpage

\acknowledgments
JRL was supported by the U.S. Department of Energy (DOE), Office of Science, Office of Workforce Development for Teachers and Scientists, Office of Science Graduate Student Research (SCGSR) program.
The SCGSR program is administered by the Oak Ridge Institute for Science and Education (ORISE) for the DOE.
ORISE is managed by ORAU under contract number DE‐SC0014664.
JRL was also supported by the National Science Foundation (NSF)
Mathematical Sciences Graduate Internship (MSGI).
GC, MRP, and DE were supported by the Earth System Model Development program area of the U.S. DOE, Office of Science, Office of Biological and Environmental Research as part of the multi-program, collaborative Integrated Coastal Modeling (ICoM) project.


\datastatement
The Python code which implements the optimization problem formulated in Section \ref{sec:optimization} is publicly available on both GitHub and Zenodo \citep{lilly2023a}.
\begin{itemize}[leftmargin=2.5cm]
    \item[GitHub:] \url{https://github.com/jeremy-lilly/rk-fb-optimization/} \\ \url{tree/32d2480f5ffcb66f4bea62653cb451b4373b6ace}
    \item[Zenodo:] \url{https://doi.org/10.5281/zenodo.7958458}
\end{itemize}
Similarly, the source code for the shallow water model used here can be found on GitHub and Zenodo \citep{engwrida2023}.
\begin{itemize}[leftmargin=2.5cm]
    \item[GitHub:] \url{https://github.com/dengwirda/swe-python/} \\ \url{tree/e7a7beb3ddcc9e27a17ce9eac2de196d77b18db8}
    \item[Zenodo:] \url{https://doi.org/10.5281/zenodo.7958483}
\end{itemize}
The data generated for this paper are also publicly available on Zenodo.
This includes the model output used to generate relative vorticity plots, model output to generate convergence plots, and log files from CFL performance experiments \citep{lilly2023b}.
\begin{itemize}[leftmargin=2.5cm]
    \item[Zenodo:] \url{https://doi.org/10.5281/zenodo.7958518}
\end{itemize}


\appendix
\renewcommand{\thesection}{A}  

\appendixtitle{Supplementary Mathematical Details}


\subsection{Local Truncation Error Analysis}
\label{subsec:lte_analysis}

What follows is a simplified local truncation error analysis of FB-RK(3,2) that shows that the method is at least \( \O\left( (\Delta t)^2 \right) \) on a particular problem for any choice of the FB-weights \( \beta_i \).
To ease notation and simplify symbolic calculations, consider a simplified version of the linearized SWEs \eqref{eqn:linswe_nzmf} with zero mean flow, in one spatial dimension, and where \( f = 0 \):
\begin{linenomath}
\begin{equation}
\begin{aligned}
    \pd{u}{t} &= -c \pd{\eta}{x} \\
    \pd{\eta}{t} &= -c \pd{u}{x} \,.
\end{aligned} \label{eqn:1d_linswe}
\end{equation}
\end{linenomath}
Assume that \( u \) is sufficiently smooth that derivatives up to sixth order in both \( t \) and \( x \) exist, and that \( \eta \) is sufficiently smooth that derivatives up to fifth order in both \( t \) and \( x \) exist. Further, assume that both \( u \) and \( \eta \) are sufficiently smooth that \( \pd{}{t}\pd[k]{}{x} u = \pd[k]{}{x}\pd{}{t} u \) and \( \pd{}{t}\pd[k]{}{x} \eta = \pd[k]{}{x}\pd{}{t} \eta \) for \( k = 1, 2\). This second assumption gives us
\begin{linenomath}
\begin{align}
    \pd[k]{u}{t} &= \begin{cases}
        -c^k \pd[k]{\eta}{x} & \text{if \( k \) is odd} \\
        c^k \pd[k]{u}{x} & \text{if \( k \) is even}
    \end{cases} \label{eqn:u_derivatives} \\
    \pd[k]{\eta}{t} &= \begin{cases}
        -c^k \pd[k]{u}{x} & \text{if \( k \) is odd} \\
        c^k \pd[k]{\eta}{x} & \text{if \( k \) is even}
    \end{cases} \label{eqn:eta_derivatives}
\end{align}
\end{linenomath}
for \( k = 1, 2 \).
Also, note that in the analysis that follows, the dependence of \( u \) and \( \eta \) on \( x \) is suppressed to simplify notation.

Now, by applying FB-RK(3,2) to \eqref{eqn:1d_linswe}, we can write
\begin{linenomath}
{\small\begin{align}
    \frac{u^{n+1} - u^n}{\Delta t} &=  - c \pd{\eta^n}{x} 
                                       + \frac{\Delta t}{2} c^2 \pd[2]{u^n}{x} 
                                       - \frac{(\Delta t)^2}{6} c^3 (\beta_3+1) \pd[3]{\eta^n}{x} 
                                       + \frac{(\Delta t)^3}{36} c^4 (-4\beta_1\beta_3 + 2\beta_1 + 9\beta_2\beta_3) \pd[4]{u^n}{x} \nonumber \\
                                       &\qquad
                                       - \frac{(\Delta t)^4}{12} c^5 \beta_2 \beta_3 \pd[5]{\eta^n}{x}
                                       + \frac{(\Delta t)^5}{36} c^6 \beta_1\beta_2\beta_3 \pd[6]{u^n}{x} \label{eqn:u_prime_est} \\
                                   &= U(u^n, \eta^n) \nonumber \\
    &\null \nonumber \\
    \frac{\eta^{n+1} - \eta^n}{\Delta t} &= - c \pd{u^n}{x}
                                            + \frac{\Delta t}{2} c^2 \pd[2]{\eta^n}{x}
                                            - \frac{(\Delta t)^2}{4} c^3 \beta_2 \pd[3]{u^n}{x}
                                            + \frac{(\Delta t)^3}{12} c^4 \beta_2 \pd[4]{\eta^n}{x}
                                            - \frac{(\Delta t)^4}{36} c^5 \beta_1 \beta_2 \pd[5]{u^n}{x} \label{eqn:eta_prime_est} \\
                                         &= H(u^n, \eta^n) \,. \nonumber
\end{align}}
\end{linenomath}
Notice that the right-hand-side of \eqref{eqn:u_prime_est} and \eqref{eqn:eta_prime_est} are equations in \( u^n \) and \( \eta^n \); call them \( U(u^n, \eta^n) \) and \( H(u^n, \eta^n) \) respectively.

The local truncation error for a given time-stepping scheme can be thought of as the error introduced by said scheme during one time-step, assuming exact data at the beginning of the time-step.
For this problem, the local truncation errors in \( u \) and \( \eta \) are defined to be the quantities \( \tau^n_u \) and \( \tau^n_\eta \) such that
\begin{linenomath}
\begin{align}
    \frac{u(t^{n+1}) - u(t^n)}{\Delta t} = U\left( u(t^n), \eta(t^n) \right) + \tau^n_u \label{eqn:lte_u_def} \\
    &\null \nonumber \\
    \frac{\eta(t^{n+1}) - \eta(t^n)}{\Delta t} = H\left( u(t^n), \eta(t^n) \right) + \tau^n_\eta \,. \label{eqn:lte_eta_def} 
\end{align}
\end{linenomath}

To calculate \( \tau^n_u \) and \( \tau^n_\eta \), rearrange \eqref{eqn:lte_u_def} and \eqref{eqn:lte_eta_def}, then expand \( u(t^{n+1}) \) and \( \eta(t^{n+1}) \) in \( t \) as Taylor polynomials centered at \( t = t^n \) to get
\begin{linenomath}
\begin{align}
    \Delta t\,\tau^n_u &= u(t^{n+1}) - u(t^n) - \Delta t \, U\left( u(t^n), \eta(t^n) \right) \nonumber \\
    &= \left( \sum_{k=0}^6 \frac{(\Delta t)^k}{k!}\pd[k]{u}{t}(t^n) + \O\left((\Delta t)^7\right) \right) - u(t^n) - \Delta t\, U\left( u(t^n), \eta(t^n) \right) \label{eqn:lte_u_full} \\
    &\null \nonumber \\
    \Delta t\,\tau^n_\eta &= \eta(t^{n+1}) - \eta(t^n) - \Delta t\, H\left( u(t^n), \eta(t^n) \right) \nonumber \\
    &= \left( \sum_{k=0}^5 \frac{(\Delta t)^k}{k!}\pd[k]{\eta}{t}(t^n) + \O\left((\Delta t)^6\right) \right) - \eta(t^n) - \Delta t\, H\left( u(t^n), \eta(t^n) \right) \,. \label{eqn:lte_eta_full}
\end{align}
\end{linenomath}

Now, simplify \eqref{eqn:lte_u_full} and \eqref{eqn:lte_eta_full}. Starting with \eqref{eqn:lte_u_full}, we get
\begin{linenomath}
{\small\begin{align}
    \Delta t\, \tau^n_u &= \left( \sum_{k=0}^3 \frac{(\Delta t)^k}{k!}\pd[k]{u}{t}(t^n) + \O\left((\Delta t)^4\right) \right) - u(t^n) - \Delta t\, U\left( u(t^n), \eta(t^n) \right) \nonumber \\
    &= \left( u(t^n)
              + \Delta t\, \underbrace{\pd{u}{t}(t^n)}_{= -c\pd{\eta}{x}(t^n)}
              + \frac{(\Delta t)^2}{2}\underbrace{\pd[2]{u}{t}(t^n)}_{= c^2 \pd[2]{u}{x}(t^n)}
              + \frac{(\Delta t)^3}{6}\pd[3]{u}{t}(t^n)
              + \O\left((\Delta t)^4\right) \right) \nonumber \\
    &\qquad + \left( - u(t^n)
                     + \Delta t\, c \pd{\eta^n}{x}(t^n)
                     - \frac{(\Delta t)^2}{2} c^2 \pd[2]{u^n}{x}(t^n)
                     + \frac{(\Delta t)^3}{6} c^3 (\beta_3+1) \pd[3]{\eta^n}{x}(t^n)
                     + \O\left((\Delta t)^4\right) \right) \nonumber \\
    &= \left( \frac{(\Delta t)^3}{6}\pd[3]{u}{t}(t^n)
              + \frac{(\Delta t)^3}{6} c^3 (\beta_3+1) \pd[3]{\eta}{x}(t^n)
              + \O\left((\Delta t)^4\right) \right) \,. \label{eqn:lte_u_no_beta_assuptions}
\end{align}}
\end{linenomath}
From \eqref{eqn:lte_u_no_beta_assuptions}, we see that \( \Delta t\,\tau^n_u = \O\left((\Delta t)^3\right) \) no matter the values of the FB-weights, so
\begin{linenomath}
\begin{equation}
    \tau^n_u = \O \left((\Delta t)^2\right). \label{eqn:lte_u_at_least_O2}
\end{equation}
\end{linenomath}

Next simplify \eqref{eqn:lte_eta_full} in a similar way:
\begin{linenomath}
{\small\begin{align}
    \Delta t\, \tau^n_\eta &= \left( \sum_{k=0}^3 \frac{(\Delta t)^k}{k!}\pd[k]{\eta}{t}(t^n) + \O\left((\Delta t)^4\right) \right) - \eta(t^n) - \Delta t H\left( u(t^n), \eta(t^n) \right) \nonumber \\
    &= \left( \eta(t^n)
              + \Delta t\, \underbrace{\pd{\eta}{t}(t^n)}_{= -c \pd{u}{x}(t^n)}
              + \frac{(\Delta t)^2}{2}\underbrace{\pd[2]{\eta}{t}(t^n)}_{= c^2\pd[2]{\eta}{x}(t^n)}
              + \frac{(\Delta t)^3}{6}\pd[3]{\eta}{t}(t^n) 
              + \O\left((\Delta t)^4\right) \right) \nonumber \\
    &\qquad + \bigg( - \eta(t^n)
                     + \Delta t\, c \pd{u^n}{x}(t^n)
                     - \frac{(\Delta t)^2}{2} c^2 \pd[2]{\eta^n}{x}(t^n)
                     + \frac{(\Delta t)^3}{4} c^3 \beta_2 \pd[3]{u^n}{x}(t^n)
                     + \O\left((\Delta t)^4\right) \bigg) \nonumber \\
    &= \left( \frac{(\Delta t)^3}{6}\pd[3]{\eta}{t}(t^n)
              + \frac{(\Delta t)^3}{4} c^3 \beta_2 \pd[3]{u}{x}(t^n)
              + \O\left((\Delta t)^4\right) \right) \,. \label{eqn:lte_eta_no_beta_assuptions}
\end{align}}
\end{linenomath}
From \eqref{eqn:lte_eta_no_beta_assuptions}, we see that with no assumptions on the FB-weights, we have \( \Delta t\tau^n_\eta = \O\left((\Delta t)^3\right) \), so
\begin{linenomath}
\begin{equation}
    \tau^n_\eta = \O\left((\Delta t)^2\right). \label{eqn:lte_eta_at_least_O2}
\end{equation}
\end{linenomath}


\subsection{Linearization of the Shallow Water Equations}
\label{subsec:linearization}

Here, we present the derivation of \eqref{eqn:nondimless_linswe_nzmf}, the linearized SWEs centered about a nonzero mean flow.
The nonlinear SWEs on a rotating sphere are given by
\begin{linenomath}
\begin{equation}
\begin{aligned}
    \pd{\breve{\u}}{t} + \left( \nabla \times \breve{\u} + f\mathbf{k} \right) \times \breve{\u} &= -\nabla\frac{\abs{\breve{\u}}^2}{2} - g\nabla \breve{h} \\
    \pd{\breve{h}}{t} + \nabla \cdot \left(\breve{h}\breve{\u}\right) &= 0 \,,
\end{aligned} \label{eqn:nonlinear_swe}
\end{equation}
\end{linenomath}
where \( \breve{\u}(x, y, t) = \left(\breve{u}(x, y, t), \breve{v}(x, y, t)\right) \) is the horizontal fluid velocity, \( x \) and \( y \) are the spatial coordinates, \( t \) is the time coordinate, \( f \) is the Coriolis parameter, \( g \) is the gravitational constant, and \( \breve{h} \) is the fluid thickness.

To linearize \eqref{eqn:nonlinear_swe} about a nonzero mean flow, assume that
\begin{linenomath}
\begin{align}
    \breve{\u}(x, y, t) &= \tilde{\U} + \tilde{\u}(x, y, t) \label{eqn:meanflow_u} \\
    \breve{h}(x, y, t) &= H + \tilde{h}(x, y, t) \,, \label{eqn:meanflow_h}
\end{align}
\end{linenomath}
where \( \tilde{\U} = (\tilde{U}, \tilde{V}) \) is the constant (with respect to time and space) mean horizontal fluid velocity, \( \tilde{\u} = \left(\tilde{u}(x, y, t), \tilde{v}(x, y, t)\right) \) is a \emph{small} perturbation of \( \tilde{\U} \), \( H \) is the constant (with respect to time and space) fluid thickness when at rest, and \( \tilde{h} = \tilde{h}(x, y, t) \) is a small perturbation of \( H \). In both cases, small means that the quantities themselves, as well as their derivatives, are negligible when multiplied together or squared.

Now, we make the substitutions \( \breve{\u} = \tilde{\U} + \tilde{\u} \) and \( \breve{h} = H + \tilde{h} \) and simplify term-by-term. The kinetic energy term becomes
\begin{linenomath}
\begin{align}
    -\nabla\frac{\abs{\breve{\u}}^2}{2} &= -\nabla\frac{1}{2} \left( \breve{u}^2 + \breve{v}^2 \right) \nonumber \\
    &= -\nabla\frac{1}{2} \left( \tilde{U}^2 + 2\tilde{U}\tilde{u} + \tilde{u}^2 + \tilde{V}^2 + 2\tilde{V}\tilde{v} + \tilde{v}^2 \right) \nonumber \\
    &= - \tilde{U}\nabla\tilde{u} - \tilde{V}\nabla\tilde{v} - \frac{1}{2}\left( \nabla\tilde{u}^2 + \nabla\tilde{v}^2 \right) \nonumber \\
    &\approx - \tilde{U}\nabla\tilde{u} - \tilde{V}\nabla\tilde{v} \nonumber \\
    &= \tilde{\U} \cdot \nabla \tilde{\u} \,, \label{eqn:simple_ke}
\end{align}
\end{linenomath}
where \( \nabla\tilde{\u} \) is the Jacobian matrix \( \begin{pmatrix} \partial_x\tilde{u} & \partial_x\tilde{v} \\ \partial_y\tilde{u} & \partial_y\tilde{v} \end{pmatrix} \) and \( \tilde{\U} \cdot \nabla \tilde{\u} \) denotes a \( 2 \times 1 \) column vector wherein the \( i \)-th entry is the dot product of \( \tilde{\U} \) with the \( i \)-th row of \( \nabla\tilde{\u} \). The potential vorticity and Coriolis term becomes
\begin{linenomath}
\begin{align}
    \left( \nabla \times \breve{\u} + f\mathbf{k} \right) \times \breve{\u} &= \left( \left(\pd{\breve{v}}{x} - \pd{\breve{u}}{y}\right)\mathbf{k} + f\mathbf{k} \right) \times \breve{\u} \nonumber \\
    &= \mathbf{i}\left( -\breve{v}\left( \pd{\breve{v}}{x} - \pd{\breve{u}}{y} + f \right) \right) + \mathbf{j}\left( \breve{u}\left( \pd{\breve{v}}{x} - \pd{\breve{u}}{y} + f \right) \right) \nonumber \\
    &= \mathbf{i}\left( -f(\tilde{V} + \tilde{v}) - \tilde{V}\left(\pd{\tilde{v}}{x} - \pd{\tilde{u}}{y}\right) - \tilde{v}\left(\pd{\tilde{v}}{x} - \pd{\tilde{u}}{y}\right) \right) \nonumber \\
    &\qquad + \mathbf{j}\left( f(\tilde{U} + \tilde{u}) + \tilde{U}\left(\pd{\tilde{v}}{x} - \pd{\tilde{u}}{y}\right) + \tilde{u}\left(\pd{\tilde{v}}{x} - \pd{\tilde{u}}{y}\right) \right) \nonumber \\
    &\approx \mathbf{i}\left( -f(\tilde{V} + \tilde{v}) - \tilde{V}\left(\pd{\tilde{v}}{x} - \pd{\tilde{u}}{y}\right) \right) + \mathbf{j}\left( f(\tilde{U} + \tilde{u}) + \tilde{U}\left(\pd{\tilde{v}}{x} - \pd{\tilde{u}}{y}\right) \right) \nonumber \\
    &= f(\tilde{\U} + \tilde{\u})^\perp + \tilde{\zeta}\tilde{\U}^\perp \,, \label{eqn:simple_pvcor}
\end{align}
\end{linenomath}
where \( \tilde{\zeta} = \pd{\tilde{v}}{x} - \pd{\tilde{u}}{y} \) is the vertical component of the vorticity \( \nabla \times \tilde{\u} \). The thickness gradient becomes
\begin{linenomath}
\begin{align}
    -g \nabla \breve{h} &= -g \nabla (H + \tilde{h}) \nonumber \\
    &= -g \nabla \tilde{h} \,. \label{eqn:simple_thickgrad}
\end{align}
\end{linenomath}
Finally, the mass advection term becomes
\begin{linenomath}
\begin{align}
    \nabla \cdot \left( \breve{h}\breve{\u} \right) &= \nabla \cdot \left( H\tilde{\U} + H\tilde{\u} + \tilde{h}\tilde{\U} + \tilde{h}\tilde{\u} \right) \nonumber \\
    &= \nabla \cdot \left( H\tilde{\u} + \tilde{h}\tilde{\U} + \tilde{h}\tilde{\u} \right) \nonumber \\
    &\approx H(\nabla \cdot \tilde{\u}) + \nabla \cdot (\tilde{h}\tilde{\U}) \,. \label{eqn:simple_massadvec}
\end{align}
\end{linenomath}
So, linearized SWEs with a non-zero mean flow are given by
\begin{linenomath}
\begin{equation}
\begin{aligned}
    &\pd{\tilde{\u}}{t} + f(\tilde{\U} + \tilde{\u})^\perp + \tilde{\zeta}\tilde{\U}^\perp = -\tilde{\U} \cdot \nabla \tilde{\u} - g \nabla \tilde{h} \\
    &\pd{\tilde{h}}{t} + H(\nabla \cdot \tilde{\u}) + \nabla \cdot (\tilde{h}\tilde{\U}) = 0 \,.
\end{aligned} \label{eqn:linswe_nzmf_appendix}
\end{equation}
\end{linenomath}


\subsection{Arakawa C-grid Discretization}
\label{subsec:c-grid_discretization}

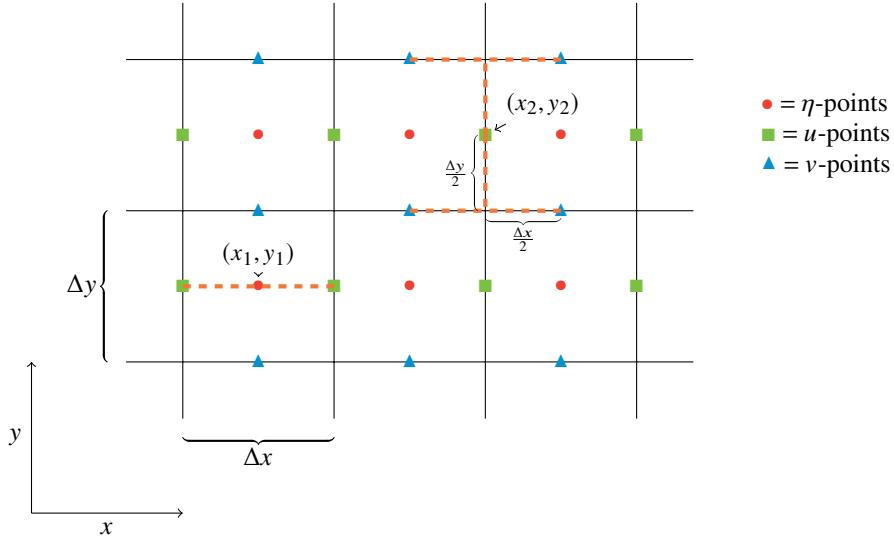
\begin{figure}[ht]
    \centering
    \begin{tikzpicture}[black]
        \def\mydot{\( \mathcolor{myred}{\bullet} \)}
        \def\mysquare{{\small \( \mathcolor{mygreen}{\blacksquare} \)}}
        \def\mytriangle{\raisebox{1px}{\small\( \mathcolor{myblue}{\blacktriangle} \)}}
    
        \draw[step=2.0] (1.25, 1.25) grid (8.75, 6.75);

        \draw[->] (0,0) -- node[midway, below]{\( x \)} ++ (2,0);
        \draw[->] (0,0) -- node[midway, left]{\( y \)} ++ (0,2);

        \foreach \x in {3, 5, 7}{
            \foreach \y in {3, 5}{
                \node at (\x, \y) {\mydot};
            }
        }

        \foreach \x in {2, 4, 6, 8}{
            \foreach \y in {3, 5}{
                \node at (\x, \y) {\mysquare};
            }
        }

        \foreach \x in {3, 5, 7}{
            \foreach \y in {2, 4, 6}{
                \node at (\x, \y) {\mytriangle};
            }
        }

        \node[align=left] (legend) at (10.5, 5) {\mydot\ \( = \eta \)-points \\ 
                                               \mysquare\ \( = u \)-points \\ 
                                               \mytriangle\ \( = v \)-points};

        \draw [thick, decorate, decoration={calligraphic brace, mirror}] (2, 1) -- node[midway, below]{\( \Delta x \)} (4,1);
        \draw [thick, decorate, decoration={calligraphic brace}] (1, 2) -- node[midway, left]{\( \Delta y \)} (1,4);

        \node (eta) at (3, 3) {};
        \node[above=1px of eta, align=center] (etalabel) {\small\( (x_1, y_1) \)};
        \draw[->] (etalabel) -- (eta);
        \draw[myorange, dashed, ultra thick] (2, 3) -- (4, 3);

        \draw[myorange, dashed, ultra thick] (6, 5) -- (6, 6);
        \draw[myorange, dashed, ultra thick] (5, 6) -- (7, 6);
        \draw[myorange, dashed, ultra thick] (6, 5) -- (6, 4);
        \draw[myorange, dashed, ultra thick] (5, 4) -- (7, 4);
        \node (u) at (6, 5) {};
        \node[above right=1px of u, align=left] (ulabel) {\small\( (x_2, y_2) \)};
        \draw[->] (ulabel) -- (u);
        \draw [decorate, decoration={calligraphic brace, mirror}] (6, 4-0.1) -- node[midway, below]{\tiny \( \frac{\Delta x}{2} \)} (7,4-0.1);
        \draw [decorate, decoration={calligraphic brace}] (6-0.1, 4) -- node[midway, left]{\tiny \( \frac{\Delta y}{2} \)} (6-0.1,5);
    \end{tikzpicture}
    \caption{An example Arakawa C-grid.}
    \label{fig:c-grid}
\end{figure}

Figure \ref{fig:c-grid} shows an example of a rectangular, staggered Arakawa C-grid \cite{arakawa1977}.
Here, the mass variable \( \eta \) is computed at cell centers and the normal components of velocity \( u \) and \( v \) are computed at the center of cell edges.

To see how to equations \eqref{eqn:dx_e} and \eqref{eqn:dy_e} are obtained, consider the following example wherein we compute a derivative in \( x \) at the point \( (x_1, y_1) \) in Figure \ref{fig:c-grid}.
On this grid, the derivative \( \pd{}{x} \) of some quantity can be approximated by taking a centered difference of values at the two nearest points:
\begin{linenomath}
\begin{align*}
    \pd{}{x} \left( e^{ikx_1+i\ell y_1} \right) &= \frac{e^{ik\left(x_1 + \frac{\Delta x}{2}\right) + i\ell y_1} - e^{ik\left(x_1 + \frac{\Delta x}{2}\right) + i\ell y_1}}{\Delta x} \\
    &= \frac{e^{ik\frac{\Delta x}{2}} - e^{-ik\frac{\Delta x}{2}}}{\Delta x}e^{ikx_1 + i\ell y_1} \\
    &= \frac{i 2 \sin\left(k\frac{\Delta x}{2}\right)}{\Delta x}e^{ikx_1 + i\ell y_1}.
\end{align*}
\end{linenomath}
An identical calculation hold for derivatives in \( y \).

To calculate the Coriolis terms, we need data for \( u \) at \( v \)-points, and for \( v \) at \( u \)-points, called \( u_v \) and \( v_u \) respectively.
To see how equations \eqref{eqn:coriolis_u} and \eqref{eqn:coriolis_v} are obtained, consider the following example wherein we calculate \( v_u \) at \( (x_2, y_2) \) in Figure \ref{fig:c-grid} by taking the average of \( v \) data at the four nearest \( v \)-points:
\begin{linenomath}
\begin{align*}
    v_u(x_2, y_2, t) &= \frac{1}{4}\hat{v}(k, \ell, t)\bigg( e^{ik\left(x_2+\frac{\Delta x}{2}\right) + i\ell\left(y_2+\frac{\Delta y}{2}\right)} 
                                                             + e^{ik\left(x_2+\frac{\Delta x}{2}\right) + i\ell\left(y_2-\frac{\Delta y}{2}\right)} \\
                                                             &\qquad\qquad\qquad\quad
                                                             + e^{ik\left(x_2-\frac{\Delta x}{2}\right) + i\ell\left(y_2+\frac{\Delta y}{2}\right)} 
                                                             + e^{ik\left(x_2-\frac{\Delta x}{2}\right) + i\ell\left(y_2-\frac{\Delta y}{2}\right)} \bigg) \\
    &= \frac{1}{4}\left( e^{ik\frac{\Delta x}{2} + i\ell\frac{\Delta y}{2}}
                         + e^{ik\frac{\Delta x}{2} - i\ell\frac{\Delta y}{2}}
                         + e^{-ik\frac{\Delta x}{2} + i\ell\frac{\Delta y}{2}}
                         + e^{-ik\frac{\Delta x}{2} - i\ell\frac{\Delta y}{2}} \right) \hat{v}(k,\ell,t)e^{ikx_2+i\ell y_2} \\
    &= \frac{1}{4} \left( 4 \cos\left(k\frac{\Delta x}{2}\right) \cos\left(\ell\frac{\Delta y}{2}\right) \right) \hat{v}(k, \ell, t)e^{ikx_2+i\ell y_2} \\
    &= \cos\left(k\frac{\Delta x}{2}\right) \cos\left(\ell\frac{\Delta y}{2}\right) \hat{v}(k, \ell, t)e^{ikx_2+i\ell y_2}.
\end{align*}
\end{linenomath}
An identical calculation holds for \( u_v \).


\subsection{Formulation of \( \tilde{G}\)}
\label{subsec:derivation_of_Gtilde}

Assume that the mean flow \( \U \) in the linearized SWEs \eqref{eqn:linswe_nzmf} is taken to be zero, and apply a Fourier transform in space.
The resulting equations can be written as
\begin{linenomath}
\begin{equation}
    \frac{\partial}{\partial t} \begin{pmatrix}
        \hat{u} \\ \hat{v} \\ \hat{\eta}
    \end{pmatrix} = \begin{pmatrix}
        0 & f & -ick \\ -f & 0 & -ic\ell \\ -ick & -ic\ell & 0
    \end{pmatrix} \begin{pmatrix}
        \hat{u} \\ \hat{v} \\ \hat{\eta}
    \end{pmatrix} \,. \label{eqn:vn_exact_linswe}
\end{equation}
\end{linenomath}
These equations are spatially continuous counterparts to \eqref{eqn:vn_linswe_nzmf}.

The general solution to this system is given by
\begin{linenomath}
\begin{equation}
    \hat{\w}(t) = e^{At}\hat{\w}(0) \,, \label{eqn:general_solution}
\end{equation}
\end{linenomath}
where \( \hat{\w}(t) = \left( \hat{u}(t), \hat{v}(t), \hat{\eta}(t) \right)^T \) and \( A \) is as given in \eqref{eqn:a}.
Now, consider the general solution of the system at  time $t = \Delta t$. Then,
\begin{linenomath}
\begin{equation}
    \hat{\w}(\Delta t) = e^{A \Delta t}\hat{\w}(0) \,. \label{eqn:exact_sol_at_dt}
\end{equation}
\end{linenomath}

Multiplying initial data \( \hat{\w}(0) = \hat{\w}^0 \) by the numerical amplification matrix \( G \) defined in \eqref{eqn:amp_matrix} returns the numerical solution after one time-step. Assume that the FB-RK(3,2) scheme produces an exact solution, that is, assume \( G\hat{\w}^0 = \hat{\w}(\Delta t) \). Then
\begin{linenomath}
\begin{equation}
    G\hat{\w}^0 = e^{A \Delta t}\hat{\w}^0 \,. \label{eqn:G_and_Gtilde}
\end{equation}
\end{linenomath}
This holds for any choice of initial data \( \hat{\w}^0 \). Take \( \hat{\w}^0 = \mathbf{e}_i \), where \( \mathbf{e}_i \) is the \(i\)-th standard basis vector for \( i = 1, 2, 3 \) to get
\begin{linenomath}
\begin{equation}
    G\mathbf{e}_i = e^{A \Delta t}\mathbf{e}_i \qquad \text{for } i = 1, 2, 3 \,. \label{eqn:cols_of_G_and_Gtilde}
\end{equation}
\end{linenomath}
Because \( G\mathbf{e}_i \) and \( e^{A \Delta t}\mathbf{e}_i \) are just the \(i\)-th columns of \( G \) and \( e^{A \Delta t} \) respectively, it follows that \( G = e^{A \Delta t} \). The reverse implication is trivial. To simplify notation, we write
\begin{linenomath}
\begin{equation}
    \tilde{G} \coloneqq e^{A \Delta t} \,. \label{eqn:Gtilde}
\end{equation}
\end{linenomath}


\subsection{Constructing balanced initial conditions}
\label{subsec:ICs}

While the shallow water cases described by \cite{galewsky2004} and \cite{williamson1992} have been successfully implemented for various structured grid models, we briefly describe an initialization procedure that is used to produce geostrophically balanced initial conditions for the unstructured TRiSK-like model used in this study.
Introducing a matrix-vector form of the nonlinear SWE system presented in \eqref{eqn:nonlinear_swe}
\begin{linenomath}
\begin{equation}
\begin{aligned}
    \pd{u}{t} + (\mathbf{C}u + f) u^\perp &= -\mathbf{G}\left( g(h + z_{b}) + K \right)\,, \\
    \pd{h}{t} + \mathbf{D}(u h) &= 0 \,, \\
    K &= \frac{1}{2}\abs{u}^2 \,,
\end{aligned} \label{eqn:nonlinear_swe_matrix}
\end{equation}
\end{linenomath}
where \( \mathbf{D} \), \( \mathbf{G} \), \( \mathbf{C} \) are sparse linear operators encoding the action of divergence, gradient and curl on a given discrete vector, as per \cite{ringler2010}.
A balanced thickness distribution \( h^0 \) can be found for a given velocity field \( u^0 \)  by setting \( \pd{u}{t} = 0 \) and taking the divergence of momentum tendencies
\begin{linenomath}
\begin{equation}
\begin{aligned}
    K^0 = \frac{1}{2}\abs{u^0}^2 \,, \\
    \mathbf{D} \big[ (\mathbf{C}u^0 + f) (u^0)^\perp + \mathbf{G} K^0\big] &= -g\mathbf{D} \mathbf{G}h^0\,.
\end{aligned} \label{eqn:ic_solve}
\end{equation}
\end{linenomath}
Recognizing the product \( \mathbf{D} \mathbf{G} = \mathbf{L} \) is an approximation of the continuous Laplacian \( \nabla^2 \), \eqref{eqn:ic_solve} can be seen as an elliptic problem that can be solved by inverting the system of linear equations
\begin{linenomath}
\begin{equation}
\begin{aligned}
    h^0 = -\frac{1}{g} \mathbf{L}^{-1} \mathbf{D} \big[ (\mathbf{C}u^0 + f) (u^0)^\perp + \mathbf{G} K^0\big] \,.
\end{aligned} \label{eqn:ic_invert}
\end{equation}
\end{linenomath}
In this study we solve (\ref{eqn:ic_invert}) using SciPy's \citep{virtanen2020} \texttt{gcrotmk} iterative solver.
Compared to the evaluation of analytical profiles, the construction of initial conditions through the solution of \eqref{eqn:ic_invert} is agnostic to details of the computational mesh used and ensures the discrete state is balanced with respect to the spatial operators used to advance the flow. This approach has proved robust for the unstructured model configurations considered in this study.



\bibliographystyle{ametsocV6}
\bibliography{references}


\end{document}